\documentclass[12pt]{amsart}
\usepackage{geometry}
\usepackage{xcolor}
\usepackage{amsaddr}
\usepackage{url}
\usepackage{graphicx}
\usepackage{float}
\usepackage{tikz}
\usepackage{cite}
\usepackage{amssymb,amsmath,amsthm}
\numberwithin{equation}{section}
\allowdisplaybreaks
\newtheorem{thm}{Theorem}[section]
\newtheorem{lemma}[thm]{Lemma}

\newtheorem{definition}[thm]{Definition}

\renewcommand{\(}{\left(}
\renewcommand{\)}{\right)}

\newcommand{\E}{{\rm E}}
\newcommand{\tr}{{\rm tr}}
\newcommand{\mb}{\mathbf}
\newcommand{\bs}{\boldsymbol}

\newcommand{\bbA}{\mathbf{A}}

\newcommand{\bbX}{\mathbf{X}}

\newcommand{\du}{\circ}

\begin{document}
\title[Test for Series correlations]{Model-Free Tests for Series Correlation in Multivariate Linear Regression}

\author{Yanqing Yin}
\address{School of Mathematics and Statistics,
Jiangsu Normal University,
Xuzhou P.R.C., 221116.}
\email{yinyq@jsnu.edu.cn}
\thanks{Yanqing Yin was partially supported by a project under Grant NSFC11801234 and Grant NSFC11701234, a Project of Natural Science Foundation of Jiangsu Province (BK20181000), a project of The Natural Science Foundation of the Jiangsu Higher Education Institutions of China (18KJB110008), the Priority Academic Program Development of Jiangsu Higher Education Institutions and a Research Support Project of Jiangsu Normal University (17XLR014).}


\maketitle
\begin{abstract}
Testing for series correlation among error terms is a basic problem in linear regression model diagnostics. The famous Durbin-Watson test and Durbin's h-test rely on certain model assumptions about the response and regressor variables. The present paper proposes simple tests for series correlation that are applicable in both fixed and random design linear regression models. The test statistics are based on the regression residuals and design matrix.  The test procedures are robust under different distributions of random errors. The asymptotic distributions of the proposed statistics are derived via a newly established joint central limit theorem for several general quadratic forms and the delta method. Good performance of the proposed tests is demonstrated by simulation results.
\end{abstract}

\section{Introduction}

Linear regression is an important topic in statistics and has been found to be useful in almost all aspects of data science, especially in business and economics statistics and biostatistics. Consider the following multivariate linear regression model
\begin{align}\label{model}
  \mathbb{Y}=\mathbb{X}'\mb \beta+{{\varepsilon}},
\end{align}
where $\mathbb{Y}$ is the response variable, $\mathbb{X}=(x_{1},x_{2},\cdots,x_{p})'$ is a $p$-dimensional vector of regressors, $\beta=\(\beta_1,\beta_2,\cdots,\beta_p\)'$ is a $p$-dimensional regression coefficient vector and $\bs\varepsilon$ is random errors with zero mean. Suppose we obtain $n$ samples from this model, that is, $\mb Y=(y_1,y_2,\cdots,y_n)'$ with design matrix $\mb X=(\mb x_1,\mb x_2,\cdots,\mb x_n)'$, where for $i=1,\cdots,n$, $\mb x_i=\(\mb x_{i1}, \mb x_{i2},\cdots,\mb x_{ip}\)'$. The first task in a regression problem is to make statistical inference about the regression coefficient vector. By applying the ordinary least squares (OLS) method, we obtain the estimate $\hat \beta=\(\mb X'\mb X\)^{-1}\mb X'\mb Y$ for coefficient vector $\beta$. In most applications of linear regression models, we need the assumption that the random errors $\{\varepsilon_i\}_{i=1}^n$ are uncorrelated and homoscedastic. That is to say, we assume $${\rm Cov}\(\varepsilon_i,\varepsilon_j\)=\begin{cases} \sigma^2& \mbox{for $i=j$,}\\
     0 &\mbox{for $i\neq j$},\end{cases}$$ where $\sigma^2$ are unknown. With this assumption, the Gauss-Markov theorem states that
the ordinary least squares estimate (OLSE) $\hat \beta$ is the best linear unbiased estimator (BLUE). When this assumption does not hold, we suffer from a loss of efficiency and, even worse, make wrong inferences in using OLS. For example, positive serial correlation in the regression error terms will typically lead to artificially small standard errors for the regression coefficient when we apply the classic linear regression method, which will cause the estimated t-statistic to be inflated, indicating significance even when there is in fact none. Therefore, tests for heteroscedasticity and series correlation are important when applying linear regression.

For detecting heteroscedasticity, in one of the most
cited papers in econometrics, White \cite{white1980heteroskedasticity} proposed a test based on  comparing the Huber-White covariance estimator to the usual covariance estimator under homoscedasticity. Many other researchers have considered this problem, for example, Breusch and Pagan \cite{breusch1979simple}, Dette and Munk \cite{Dette1998a}, Glejser \cite{glejser1969new}, Harrison and McCabe \cite{harrison1979test}, Cook
and Weisberg \cite{cook1983diagnostics}, and Azzalini and Bowman\cite{azzalini1993use}. Recently, Li and Yao \cite{li2015homoscedasticity} and Bai, Pan and Yin \cite{yin2018test} proposed tests for heteroscedasticity that are valid in both low- and high-dimensional regressions. Their tests were shown by simulations to perform better than some classic tests.

The most famous test for series correlation, the Durbin-Watson test, was proposed in \cite{durbin1950testing,durbin1951testing,durbin1971testing}. The Durbin-Watson test statistic is based on the residuals $e_1,e_2,\cdots,e_n$ from linear regression. The researchers considered the statistic
$$d=\frac{\sum_{i=2}^n\(e_i-e_{i-1}\)^2}{\sum_{i=1}^n e_i^2},$$ whose small-sample distribution was derived by John von Neumann. In the original papers,  Durbin and Watson investigated the distribution of this statistic under the classic independent framework, described the test procedures and provided tables of the
bounds of significance.
 However, the asymptotic results were derived under the normality assumption on the error term, and as noted by Nerlove and Wallis \cite{Nerlove1966Use}, although the Durbin-Watson test appeared to work well in an independent observations framework, it may be asymptotically
biased and lead to inadequate conclusions for linear regression models containing
lagged dependent random variables. New alternative test procedures, for instance, Durbin's h-test and t-test \cite{Durbin1970Testing}, were proposed to address this problem; see also Inder \cite{Inder1986An}, King and Wu \cite{king1991small}, Stocker \cite{Stocker2007On}, Bercu and Pro\"{\i}a \cite{Bercu2013A}, Gen\c{c}ay and Signori
\cite{Gen2015Multi} and Li and Gen\c{c}ay \cite{Li2017Tests} and references therein. However, all these tests were proposed under some model assumptions on the regressors and/or the response variable. Moreover, Durbin's h-test requires a Gaussian distribution of the error term.  Thus, some common models are excluded. In fact, since it is difficult to assess whether the regressors and/or the response are lag dependent, model-free tests for the regressors and response variable appear to be appropriate.

The present paper proposes a simple test procedure without assumptions on the response variable and regressors that is valid in both low- and high-dimensional multivariate linear regression. The main idea, which is simple but proves to be useful, is to express the mean and variance of the test statistic by making use of the residual maker matrix.  In addition to a general joint central limit theorem for several quadratic forms, which is proved in this paper and may have its own interest, we consider a Box-Pierce-type test for series correlation. Monte Carlo simulations show that our test procedures perform well in situations where some classic test procedures are inapplicable.

\section{Test for series correlation in linear regression model}

\subsection{Notation}

Let
$\bbX$ be the design matrix,
and let $\mb R=(r_{i,j})=\mb I_n-\mb X(\mb X'\mb X)^{-1}\mb X'$ be the residual maker matrix, where $\mb H=(h_{i,j})=\mb X(\mb X'\mb X)^{-1}\mb X'$ is the
hat matrix (also known as the projection matrix).
We assume that the noise vector $\bs \varepsilon=\mb \Sigma^{1/2}\bs \epsilon,$ where $\bs \epsilon$ is an $n$-dimensional random vector whose entries $\epsilon_1,\cdots,\epsilon_n$ are independent with zero means, unit variances and the same finite fourth-order moments $M_4$, and $\mb \Sigma^{1/2}$ is an n-dimensional nonnegative definite nonrandom matrix with bounded spectral norm.
Then, the OLS residuals are ${\mb e}=(e_1,\cdots,e_n)'=\mb R\mb \Sigma^{1/2} \bs\epsilon$. We note that we will use $\du$ to indicate the Hadamard product of two matrices in the rest of this paper.

\subsection{Test for a given order series correlation}\label{dwtest}
To test for a given order series correlation, for any number $q\leq n$, denote
$$\gamma_{\tau}=\sum_{i={\tau+1}}^ne_ie_{i-\tau}=\mb e'\mb P_{\tau}\mb e=\bs \epsilon'\mb\Sigma^{1/2}\mb R\mb P_{\tau}\mb R\mb\Sigma^{1/2}\bs \epsilon, \quad \tau=0,\cdots,q,$$
where $\mb P_{\tau}=\(p_{i,j}^{(\tau)}\)_{n\times n}$ with $p_{i,j}^{(\tau)}=\begin{cases} 1 & \mbox{for $i-j=\tau$,}\\
     0 &\mbox{other}.\end{cases}$

First, for $0\leq \tau\leq  q,$ we have
\begin{align}\label{mean}
\E \gamma_{\tau}=\E\tr\bs \epsilon'\mb\Sigma^{1/2}\mb R\mb P_{\tau}\mb R\mb\Sigma^{1/2}\bs \epsilon=\tr\mb\Sigma^{1/2}\mb R\mb P_{\tau}\mb R\mb\Sigma^{1/2}=\tr\mb R\mb P_{\tau}\mb R\mb\Sigma.
\end{align}
Denote $\mb\Sigma^{1/2}\mb R\mb P_{\tau_1}\mb R\mb\Sigma^{1/2}=(a_{i,j})$ and $\mb\Sigma^{1/2}\mb R\mb P_{\tau_2}\mb R\mb\Sigma^{1/2}=(b_{i,j})$, and set $\nu_4=M_4-3;$ we then have, for $\tau_1,\tau_2=0,\cdots,q,$
\begin{align}\label{cov}
&{\rm Cov}(\gamma_{\tau_1},\gamma_{\tau_2})\\\notag
=&\E\(\bs \epsilon'\mb\Sigma^{1/2}\mb R\mb P_{\tau_1}\mb R\mb\Sigma^{1/2}\bs \epsilon-\tr\mb R\mb P_{\tau_1}\mb R\mb\Sigma\)\(\bs \epsilon'\mb\Sigma^{1/2}\mb R\mb P_{\tau_2}\mb R\mb\Sigma^{1/2}\bs \epsilon-\tr\mb R\mb P_{\tau_2}\mb R\mb\Sigma\)\\\notag
=&\E\(\sum_{i}a_{i,i}\(\epsilon_i^2-1\)+\sum_{i\neq j}a_{i,j}\epsilon_i\epsilon_j\)\(\sum_{i}b_{i,i}\(\epsilon_i^2-1\)+\sum_{i\neq j}b_{i,j}\epsilon_i\epsilon_j\)\\\notag
=&\nu_4\tr\(\(\mb\Sigma^{1/2}\mb R\mb P_{\tau_1}\mb R\mb\Sigma^{1/2}\)\du\(\mb\Sigma^{1/2}\mb R\mb P_{\tau_2}\mb R\mb\Sigma^{1/2}\)\)+\tr\(\mb\Sigma^{1/2}\mb R\mb P_{\tau_1}\mb R\mb\Sigma^{1/2}\)\(\mb\Sigma^{1/2}\mb R\mb P_{\tau_2}\mb R\mb\Sigma^{1/2}\)\\\notag
&+\tr\(\mb\Sigma^{1/2}\mb R\mb P_{\tau_1}\mb R\mb\Sigma^{1/2}\)\(\mb\Sigma^{1/2}\mb R\mb P_{\tau_2}\mb R\mb\Sigma^{1/2}\)'\\\notag
=&\nu_4\tr\(\(\mb\Sigma^{1/2}\mb R\mb P_{\tau_1}\mb R\mb\Sigma^{1/2}\)\du\(\mb\Sigma^{1/2}\mb R\mb P_{\tau_2}\mb R\mb\Sigma^{1/2}\)\)+\tr\mb P_{\tau_1}\mb R\mb\Sigma\mb R\mb P_{\tau_2}\mb R\mb\Sigma\mb R+\tr\mb P_{\tau_1}\mb R\mb\Sigma\mb R\mb P'_{\tau_2}\mb R\mb\Sigma\mb R.
\end{align}

Note that ${\rm Cov}\(\bs \varepsilon,\bs \varepsilon\)=\mb \Sigma.$
We want to test the hypothesis for $1\leq \tau\leq q,$
$$H_0: \mb\Sigma=\sigma^2\mb I_n, \quad {\mbox {where} \ } 0<\sigma^2<\infty, $$ against $$H_{1,\tau}: {\rm Cov}(\varepsilon_i,\varepsilon_{i-\tau})=\rho\neq 0.$$

Under the null hypothesis, due to (\ref{mean}) and (\ref{cov}), we obtain
\begin{align}
\E \gamma_{\tau}=\sigma^2\tr\mb P_{\tau}\mb R,
\end{align}
and
\begin{align}
  &{\rm Cov}(\gamma_{\tau_1},\gamma_{\tau_2})
  =\sigma^4\(\nu_4\tr\(\(\mb R\mb P_{\tau_1}\mb R\)\du\(\mb R\mb P_{\tau_2}\mb R\)\)+\tr\mb P_{\tau_1}\mb R\mb P_{\tau_2}\mb R+\tr\mb P_{\tau_1}\mb R\mb P'_{\tau_2}\mb R\).
\end{align}
Specifically, we have $\E \gamma_{\tau_0}=\sigma^2\tr\mb R=\sigma^2(n-p)$ and \begin{align}
  &{\rm Var}(\gamma_{\tau_0})
  =\sigma^4\(\nu_4\tr\(\mb R\du\mb R\)+2(n-p)\).
\end{align}
The validity of our test procedure requires the following mild assumptions.
\begin{description}
  \item[(1): Assumption on $p$ and $n$] The number of regressors $p$ and the sample size $n$ satisfy that $p/n\to c\in [0,1)$ as $n\to\infty$.
  \item[(2): Assumption on errors] The fourth-order cumulant of the error distribution $\nu_4\neq -2$.
\end{description}
Assumption $(2)$ excludes the rare case where the random errors are drawn from a two-point distribution with the same masses $1/2$ at $--1$ and $1$. However, if this situation occurs, our test remains valid if the design matrix satisfies the mild condition that $$\limsup_{n\to \infty}\frac{\tr \mb R\du\mb R}{n-p}=\limsup_{n\to \infty}\frac{\sum_{i=1}^nr_{i,i}^2}{\sum_{i=1}^nr_{i,i}}< 1.$$
These assumptions ensure that ${\rm Var}(\gamma_{\tau_0})$ has the same order as $n$ as $n\to \infty$, thus satisfying the condition assumed in Theorem \ref{th:CLT}.

Define $$m_{\tau}=\frac{\E\gamma_{\tau}}{\sigma^2}=\tr\mb P_{\tau}\mb R,\quad v_{\tau_1\tau_2}=\frac{{\rm Cov}(\gamma_{\tau_1},\gamma_{\tau_2})}{n\sigma^4}.$$

By applying Theorem \ref{th:CLT} presented  in Section \ref{clt}, we obtain that for $1\leq \tau\leq q,$

\begin{align*}
\frac{1}{\sqrt{n}}\(\left(
          \begin{array}{c}
            \gamma_{\tau} \\
             \gamma_0\\
          \end{array}
        \right)-\left(
                  \begin{array}{c}
                    m_{\tau} \\
                    m_0 \\
                  \end{array}
                \right)
        \)\sim N\(\mb 0, \left(
                          \begin{array}{cc}
                            v_{{\tau}{\tau}} & v_{{\tau}0} \\
                            v_{0{\tau}} & v_{00} \\
                          \end{array}
                        \right)
        \).
\end{align*}
Then, by the delta method, we obtain, as $n\to\infty$,
\begin{align*}
T_{\tau}=\frac{\sqrt{n}\(\frac{\gamma_{\tau}}{\gamma_0}-\mu_{{\tau}0}\)}{\sigma_{{\tau}0}}\sim N(0,1),
\end{align*}
where $\mu_{{\tau}0}=\frac{m_{\tau}}{m_0}$ and
\begin{align}
  \sigma_{{\tau}0}^2=&\left(
                \begin{array}{cc}
                  \frac{1}{m_{0}/n} & \frac{m_{1}/n}{\(m_{0}/n\)^2} \\
                \end{array}
              \right)
  \left(
                \begin{array}{cc}
v_{{1}{1}} & v_{{1}0} \\
                            v_{0{1}} & v_{00} \\
                \end{array}
              \right)\left(
                       \begin{array}{c}
                         \frac{1}{m_{0}/n} \\
                         \frac{m_{1}/n}{\(m_{0}/n\)^2} \\
                       \end{array}
                     \right)\\\notag
                     =&n^2\left(
                \begin{array}{cc}
                  \frac{1}{n-p} & \frac{m_{1}}{(n-p)^2} \\
                \end{array}
              \right)
  \left(
                \begin{array}{cc}
v_{{1}{1}} & v_{{1}0} \\
                            v_{0{1}} & v_{00} \\
                \end{array}
              \right)\left(
                       \begin{array}{c}
                         \frac{1}{n-p} \\
                         \frac{m_{1}}{(n-p)^2} \\
                       \end{array}
                     \right).
\end{align}
We reject $H_0$ in favor of $H_{1,\tau}$ if a large $|T_\tau|$ is observed.

\subsection{A portmanteau test for series correlation}\label{bptest}
In time series analysis, the Box-Pierce test proposed in \cite {box1970distribution} and the Ljung-Box statistic proposed in \cite{ljung1978measure} are two portmanteau tests of whether any of a group of autocorrelations of a time series are different from zero. For a linear regression model, consider the following hypothesis
$$H_0: \mb\Sigma=\sigma^2\mb I_n,$$ against $$H_{1}: {\mbox {there\ exist}}\  1\leq \tau\leq q \ {\mbox {such\ that}}\ {\rm Cov}(\varepsilon_i,\varepsilon_{i-\tau})=\rho\neq 0.$$
Applying Theorem \ref{th:CLT} and the delta method, we shall now consider the following asymptotically standard normally distributed statistic
\begin{align*}
T(q)=\frac{\sqrt{n}\(\sum_{\tau=1}^q\({2-\frac{\gamma_{\tau}}{\gamma_0}}\)^2-\sum_{\tau=1}^q(2-{\mu_{\tau 0}})^2\)}{\sigma_T}\sim N(0,1),
\end{align*}
as $n\to \infty,$
where $\mu_{\tau 0}=\frac{m_{\tau}}{m_0}$ and
$\sigma_T=\sqrt{{\nabla}' {\Sigma}_T\nabla}$ with $${\nabla}=\(\frac{2n\sum_{\tau=1}^qm_\tau (2-\frac{m_{\tau}}{n-p})}{(n-p)^2}, \frac{-2n(2-\frac{m_{1}}{n-p})}{(n-p)},\frac{-2n(2-\frac{m_{2}}{n-p})}{(n-p)},\cdots,\frac{-2n(2-\frac{m_{q}}{n-p})}{(n-p)}\)',$$
and
$${\Sigma}_T=\(v_{i,j}\)_{(q+1)\times(q+1)}, \quad i,j=0,\cdots,q.$$
Then, we reject $H_0$ in favor of $H_{1}$ if $|T(q)|$ is large.
\subsection{Discussion of the statistics} In the present subsection, we discuss the asymptotic parameters of the two proposed statistics.

If the entries in design matrix $\mb X$ are assumed to be i.i.d. standard normal, then we know that as $n\to \infty$, the diagonal entries in the  symmetric and idempotent matrices $\mb H$ and $\mb R=\mb I_n-\mb H$ are of constant order while the off-diagonal entries are of order $n^{-1/2}$. Then, the order of $m_\tau=\tr\mb P_{\tau}\mb R$ for a given $\tau>0$ is at most $n^{1/2}$ since it is exactly the summation of the $n-\tau$ off-diagonal entries of $\mb R$. Thus, elementary analysis shows that $\sigma_{{\tau}0}^2=\frac{n^2v_{11}}{(n-p)^2}+o(1)$.

For a fixed design or a more general random design, it become almost impossible to study matrices $\mb H$ and $\mb R=\mb I_n-\mb H$, except for some of the elementary properties. Thus, for the purpose of obtaining an accurate statistical inference, we suggest the use of the original parameters since we have little information on the distribution of the regressors in a fixed design, and the calculation of those parameters is not excessively complex.
\section{Simulation studies}
In this section, Monte Carlo simulations are conducted to investigate the performance of our proposed tests.
\subsection{Performance of test for first-order series correlation}
First, we consider the test for first-order series correlation of the error terms in multivariate linear regression model (\ref{model}). Note that although our theory results were derived by treating the design matrix as a constant matrix, we also need to obtain a design matrix under a certain random model in the simulations. We thus consider the situation where the regressors $\mb x_1,\mb x_2,\cdots,\mb x_f$ are lagged dependent. Formally, for a given $f$, we set $$x_{t,j}=r x_{t-1,j}+u_t, \quad j=1,\cdots,f, \quad t=1,\cdots,n,$$
where $r=0.2$ and $\{u_t\}$ are independently drawn from N(0,1). While $\{x_{i,j}\}, f+1\leq j\leq p, 1\leq i\leq n$ are independently chosen from a Student's t-distribution with 5 degrees of freedom.
The random errors $\varepsilon$ obey (1) the normal distribution N(0,1) and (2) the uniform distribution U(-1,1).
The significant level is set to $\alpha=0.05.$ Table \ref{table1} and Table \ref{table2} show the empirical size of our test (denoted as $``{\rm FDWT}"$) for different $p,n,f$ under the two error distributions.
To investigate the power of our test, we randomly choose a $\varepsilon_0$ and consider the following  AR(1) model:
$$\varepsilon_t=\rho\varepsilon_{t-1}+\vartheta_t, \quad t=1,\cdots,n$$
where $\{\vartheta_t\}$ are  independently drawn from (1) N(0,1) and (2) U(-1,1).
Tables \ref{table3} and \ref{table4} show the empirical power of our proposed test for different $p,n,f,\rho$ under the two error distributions.

These simulation results show that our test always has good size and power when $n-p$ is large and is thus applicable under the framework that $p/n\to [0,1)$ as $n\to \infty$.

\begin{table}[!hbp]
  \center
    \begin{tabular}{|c|c|c|c|c|c|}
       \hline
       $p,n$ &f   & FDWT &  $p,n$ &f & FDWT  \\
       \hline
       2,32 & 1  & 0.0486 & 8,32 & 2  &0.0428\\
       \hline
       8,32 & 4  & 0.0410 & 8,32 &8  & 0.0434\\
       \hline
       16,64 & 4 & 0.0446 & 16,64 & 12   & 0.0463\\
       \hline
       32,64 & 12   & 0.0420 & 32,64 & 24   & 0.0414\\
       \hline
       32,128 &12  & 0.0470 & 32,128 & 24 & 0.0478\\
       \hline
       64,128 &12  & 0.0479 & 64,128 & 36   & 0.0430\\
       \hline
       128,256 & 12 & 0.0509 & 128,256 & 24  &0.0486\\
       \hline
       128,256 & 64  & 0.0504 & 128,256 & 128  &0.0422\\
       \hline
       128,512 & 24 & 0.0519 & 128,512 & 64  &0.0496\\
       \hline
       128,512 & 96 & 0.0487 & 128,512 & 128  &0.0497\\
       \hline
       256,512 & 64  & 0.0469 & 256,512 & 96  &0.0492\\
       \hline
       256,512 & 144  & 0.0472 & 256,512 & 256  &0.0486\\
       \hline
       256,1028 & 64  & 0.0457 & 256,1028 & 96  &0.0498\\
       \hline
       256,1028 & 144  & 0.0473 & 256,1028 & 256  &0.0487\\
       \hline
       512,1028 & 12  & 0.0463 & 512,1028 & 96 &0.0506\\
       \hline
       512,1028 & 144  & 0.0520 & 512,1028 & 256  &0.0478\\
       \hline
       512,1028 & 288  & 0.0460 & 512,1028 & 314  &0.0442\\
       \hline
       512,1028 & 440  & 0.0438 & 512,1028 & 512  &0.0443\\
       \hline
     \end{tabular}
  \caption{Empirical size under Gaussian error assumption}\label{table1}
  \end{table}

\begin{table}[!hbp]
  \center
    \begin{tabular}{|c|c|c|c|c|c|}
       \hline
       $p,n$ &f & FDWT &  $p,n$ &f  & FDWT  \\
       \hline
       2,32 & 1  & 0.0410 & 2,32 & 2 &0.0421\\
       \hline
       8,32 & 4 & 0.0414 & 8,32 &8  & 0.0468\\
       \hline
       16,64 & 4 & 0.0467 & 16,64 & 12  & 0.0450\\
       \hline
       32,64 & 12  & 0.0450 & 32,64 & 24 & 0.0419\\
       \hline
       32,128 &12  & 0.0456 & 32,128 & 24  & 0.0458\\
       \hline
       64,128 &12   & 0.0479 & 64,128 & 36   & 0.0460\\
       \hline
       128,256 & 12 & 0.0509 & 128,256 & 24  &0.0476\\
       \hline
       128,256 & 64 & 0.0461 & 128,256 & 128 &0.0412\\
       \hline
       128,512 & 24 & 0.0497 & 128,512 & 64 &0.0505\\
       \hline
       128,512 & 96  & 0.0508 & 128,512 & 128  &0.0501\\
       \hline
       256,512 & 64 & 0.0525 & 256,512 & 96&0.0455\\
       \hline
       256,512 & 144  & 0.0443 & 256,512 & 256 &0.0461\\
       \hline
       256,1028 & 64  & 0.0509 & 256,1028 & 96 &0.0455\\
       \hline
       256,1028 & 144 & 0.0482 & 256,1028 & 256  &0.0465\\
       \hline
       512,1028 & 12  & 0.0491 & 512,1028 & 96 &0.0461\\
       \hline
       512,1028 & 144 & 0.0483 & 512,1028 & 256  &0.0480\\
       \hline
       512,1028 & 288  & 0.0447 & 512,1028 & 314  &0.0468\\
       \hline
       512,1028 & 440  & 0.0453 & 512,1028 & 512  &0.0459\\
       \hline
     \end{tabular}
  \caption{Empirical size under uniform distribution U(-1,1) error assumption}\label{table2}
  \end{table}

\begin{table}[!hbp]
  \center
    \begin{tabular}{|c|c|ccc|c|c|ccc|}
       \hline
       $p,n$ &f & $\rho=0.2$ & $\rho=-0.3$ & $\rho=0.5$ &  $p,n$ &f & $\rho=0.2$ & $\rho=-0.3$ & $\rho=0.5$  \\
       \hline
       2,32 & 1 & 0.1363&0.2550 & 0.5409 & 8,32 & 2 & 0.1056 & 0.1705&0.3224\\
       \hline
       8,32 & 4 & 0.0906 & 0.1724 & 0.3841 & 8,32 &8 & 0.1093 & 0.1831 & 0.3597\\
       \hline
       16,64 & 4 & 0.1888 & 0.3672 & 0.6783 & 16,64 & 12 & 0.1987 & 0.3764& 0.7199\\
       \hline
       32,64 & 12 & 0.1055 & 0.1584 & 0.3542 & 32,64 & 24 & 0.1030 & 0.1739 & 0.3791\\
       \hline
       32,128 &12 & 0.3673 & 0.6637 & 0.9552 & 32,128 & 24 & 0.3706 & 0.6655 & 0.9556\\
       \hline
       64,128 &12 & 0.1754 & 0.3335 & 0.6345 & 64,128 & 36 & 0.1897 & 0.3519 & 0.6639\\
       \hline
       128,256 & 12 & 0.3255&0.6104 & 0.9160 & 128,256 & 24 & 0.3324 & 0.6037&0.9225\\
       \hline
       128,256 & 64 & 0.3362&0.6200 & 0.9345 & 128,256 & 128 & 0.3362 & 0.6515&0.9438\\
       \hline
       128,512 & 24 & 0.9064&0.9981 & 1.0000 & 128,512 & 64 & 0.9151 & 0.9976&1.0000\\
       \hline
       128,512 & 96 & 0.9167&0.9981 & 1.0000 & 128,512 & 128 & 0.9196 & 0.9981&1.0000\\
       \hline
       256,512 & 64 & 0.5880&0.8951 & 0.9975 & 256,512 & 96 & 0.6041 & 0.9029&0.9980\\
       \hline
       256,512 & 144 & 0.6019&0.8963 & 0.9990 & 256,512 & 256 & 0.6117 & 0.9103&0.9987\\
       \hline
       256,1028 & 64 & 0.9970&1.0000 & 1.0000 & 256,1028 & 96 & 0.9973 & 1.0000&1.0000\\
       \hline
       256,1028 & 144 & 0.9971&1.0000 & 1.0000 & 256,1028 & 256 & 0.9976 & 1.0000&1.0000\\
       \hline
       512,1028 & 12 & 0.8766&0.9957 & 1.0000 & 512,1028 & 96 & 0.8829 & 0.9958&1.0000\\
       \hline
       512,1028 & 144 & 0.9201&0.9979 & 1.0000 & 512,1028 & 256 & 0.8967 & 0.9954&1.0000\\
       \hline
       512,1028 & 288 & 0.9125&0.9986 & 1.0000 & 512,1028 & 314 & 0.8946 & 0.9969&1.0000\\
       \hline
       512,1028 & 440 & 0.8942&0.9975 & 1.0000 & 512,1028 & 512 & 0.8937 & 0.9979&1.0000\\
       \hline
     \end{tabular}
  \caption{Empirical power under Gaussian error assumption}\label{table3}
  \end{table}

\begin{table}[!hbp]
  \center
        \begin{tabular}{|c|c|ccc|c|c|ccc|}
       \hline
       $p,n$ &f & $\rho=0.2$ & $\rho=-0.3$ & $\rho=0.5$ &  $p,n$ &f & $\rho=0.2$ & $\rho=-0.3$ & $\rho=0.5$  \\
       \hline
       2,32 & 1 & 0.1457&0.2521 & 0.5548 & 8,32 & 2 & 0.1245 & 0.1721&0.3478\\
       \hline
       8,32 & 4 & 0.1245 & 0.1754 & 0.3548 & 8,32 &8 & 0.1254 & 0.1845 & 0.3547\\
       \hline
       16,64 & 4 & 0.1987 & 0.3789 & 0.6567 & 16,64 & 12 & 0.1879 & 0.3478& 0.7456\\
       \hline
       32,64 & 12 & 0.1145 & 0.1544 & 0.3582 & 32,64 & 24 & 0.1125 & 0.1555 & 0.3548\\
       \hline
       32,128 &12 & 0.3825 & 0.6647 & 0.9845 & 32,128 & 24 & 0.3845 & 0.6789 & 0.9677\\
       \hline
       64,128 &12 & 0.1863 & 0.3765 & 0.6748 & 64,128 & 36 & 0.1758 & 0.3877 & 0.6478\\
       \hline
       128,256 & 12 & 0.3358&0.5978 & 0.9185 & 128,256 & 24 & 0.3495 & 0.6657&0.9244\\
       \hline
       128,256 & 64 & 0.3378&0.5899 & 0.9578 & 128,256 & 128 & 0.3392 & 0.6788&0.9584\\
       \hline
       128,512 & 24 & 0.9114&0.9945 & 1.0000 & 128,512 & 64 & 0.9121 & 0.9944&1.0000\\
       \hline
       128,512 & 96 & 0.9102&0.9977 & 1.0000 & 128,512 & 128 & 0.9157 & 0.9945&0.9999\\
       \hline
       256,512 & 64 & 0.6053&0.8979 & 0.9969 & 256,512 & 96 & 0.6020 & 0.9456&0.9978\\
       \hline
       256,512 & 144 & 0.6151&0.8966 & 1.0000 & 256,512 & 256 & 0.6135& 0.9678&1.0000\\
       \hline
       256,1028 & 64 & 0.9975&1.0000 & 1.0000 & 256,1028 & 96 & 0.9972 & 1.0000&1.0000\\
       \hline
       256,1028 & 144 & 0.9921&1.0000 & 1.0000 & 256,1028 & 256 & 0.9982 & 1.0000&1.0000\\
       \hline
       512,1028 & 12 & 0.8787&0.9944 & 1.0000 & 512,1028 & 96 & 0.8800 & 0.9976&1.0000\\
       \hline
       512,1028 & 144 & 0.9201&0.9913 & 1.0000 & 512,1028 & 256 & 0.8881 & 0.9964&1.0000\\
       \hline
       512,1028 & 288 & 0.9165&0.9959 & 1.0000 & 512,1028 & 314 & 0.8957 & 0.9967&1.0000\\
       \hline
       512,1028 & 440 & 0.8978&0.9944 & 1.0000 & 512,1028 & 512 & 0.8959 & 0.9947&1.0000\\
       \hline
     \end{tabular}
  \caption{Empirical power under uniform distribution U(-1,1) error assumption}\label{table4}
  \end{table}

\subsection{Performance of the Box-Pierce type test}
This subsection investigates the performance of our proposed Box-Pierce type test statistic $T(q)$ in subsection \ref{bptest}. The design matrix $\mb X$ is obtained in the same way as in the last subsection, with $f=p/2$, and the random error terms are assumed to obey a (1) normal distribution N(0,1) and a (2) gamma distribution with parameters 4 and 1/2. Table \ref{table5} and Table \ref{table6} show the empirical size of our test with different $n,p,q$ under the two error distributions. We consider the following AR(2) model to assess the power:
$$\varepsilon_t=\rho_1\varepsilon_{t-1}+\rho_2\varepsilon_{t-2}+\vartheta_t, \quad t=1,\cdots,n$$
where $\{\vartheta_t\}$ are  independently drawn from (1) N(0,1) and (2) Gamma(4,1/2). The design matrix $\mb X$ is obtained in the same way as before, with $f=p/2$.
Tables \ref{table7} and \ref{table8} show the empirical power of our proposed test for different $p,n,\rho$ under the two error distributions.

As shown by these simulation results, the empirical size and empirical power of the portmanteau test improve as $n-p$ tends to infinity.
\begin{table}[!hbp]
  \center
    \begin{tabular}{|c|c|cc|c|c|cc|}
       \hline
       $p,n$  & $n-p$& $q=3$  &$q=5$  &  $p,n$  & $n-p$  & $q=3$   & $q=5$  \\
       \hline
       2,32  & 30 & 0.0389 & 0.0402  & 8,32 & 24 & 0.0351 & 0.0350\\
       \hline
       16,32 & 16 & 0.0299 &0.0349 & 24,32 & 8 & 0.0208 & 0.0132\\
       \hline
       2,64 & 62 & 0.0443  &0.0505 & 32,64 & 32& 0.0391 & 0.0420\\
       \hline
       32,128 &96 & 0.0436 &0.0501 & 64,128& 64 & 0.0402 & 0.0427\\
       \hline
       32,256 & 224 & 0.0489&0.0470 & 64,256 & 192 & 0.0475 & 0.0485\\
       \hline
       128,256 & 128 & 0.0452&0.0477& 16,512 & 496 & 0.0499 & 0.0494\\
       \hline
       64,512 & 448 & 0.0490&0.0486&  128,512 & 384 & 0.0502 & 0.0513\\
       \hline
       256,512 & 256 & 0.0473&0.0438 & 64,1028 & 964 & 0.0461 & 0.0494\\
       \hline
       128,1028 & 900 & 0.0480&0.0485 & 256,1028 &772 & 0.0492 & 0.0501\\
       \hline
     \end{tabular}
  \caption{Empirical size under Gaussian error assumption}\label{table5}
  \end{table}

\begin{table}[!hbp]
  \center
    \begin{tabular}{|c|c|cc|c|c|cc|}
       \hline
       $p,n$  & $n-p$& $q=3$  &$q=5$  &  $p,n$  & $n-p$  & $q=3$   & $q=5$  \\
       \hline
       2,32  & 30 & 0.0359 & 0.0374  & 8,32 & 24 & 0.0390 & 0.0383\\
       \hline
       16,32 & 16 & 0.0265 &0.0281 & 24,32 & 8 & 0.0129 & 0.0087\\
       \hline
       2,64 & 62 & 0.0444  &0.0426 & 32,64 & 32& 0.0385 & 0.0365\\
       \hline
       32,128 &96 & 0.0430 &0.0448 & 64,128& 64 & 0.0439 & 0.0417\\
       \hline
       32,256 & 224 & 0.0497&0.0437 & 64,256 & 192 & 0.0509 & 0.0514\\
       \hline
       128,256 & 128 & 0.0487&0.0465& 16,512 & 496 & 0.0504 & 0.0498\\
       \hline
       64,512 & 448 & 0.0479&0.0511&  128,512 & 384 & 0.0498 & 0.0458\\
       \hline
       256,512 & 256 & 0.0518&0.0523 & 64,1028 & 964 & 0.0500 & 0.0489\\
       \hline
       128,1028 & 900 & 0.0490&0.0513 & 256,1028 &772 & 0.0439 & 0.0503\\
       \hline
     \end{tabular}
  \caption{Empirical size under Gamma(4,1/2) error assumption}\label{table6}
  \end{table}

\begin{table}[!hbp]
  \center
    \begin{tabular}{|c|c|cc|c|c|cc|}
       \hline
         & & $\rho_1=0.2$  &$\rho_1=0$  &    &  & $\rho_1=0.2$   & $\rho_1=0$  \\

        $p,n$  & $n-p$& $\rho_2=-0.3$  &$\rho_2=0.3$  &  $p,n$  & $n-p$  & $\rho_2=-0.3$   & $\rho_2=0.3$  \\
       \hline
       2,32  & 30 & 0.2630 & 0.1960  & 8,32 & 24 & 0.1699 & 0.1265\\
       \hline
       16,32 & 16 & 0.0890 &0.0694 & 24,32 & 8 & 0.0760 & 0.0205\\
       \hline
       2,64 & 62 & 0.5698  &0.4064 & 32,64 & 32& 0.1708 & 0.1210\\
       \hline
       32,128 &96 & 0.6660 &0.4775 & 64,128& 64 & 0.2764 & 0.2232\\
       \hline
       32,256 & 224 & 0.9849&0.9278 & 64,256 & 192 & 0.9369& 0.8167\\
       \hline
       128,256 & 128 & 0.6147&0.4335& 16,512 & 496 & 1.0000 & 1.0000\\
       \hline
       64,512 & 448 & 1.0000&1.0000&  128,512 & 384 & 0.9991 & 0.9897\\
       \hline
       256,512 & 256 & 0.9155&0.7551 & 64,1028 & 964 & 1.0000 & 1.0000\\
       \hline
       128,1028 & 900 & 1.0000&1.0000 & 256,1028 &772 & 1.0000 & 1.0000\\
       \hline
     \end{tabular}
  \caption{Empirical power under Gaussian error assumption}\label{table7}
  \end{table}

\begin{table}[!hbp]
  \center
    \begin{tabular}{|c|c|cc|c|c|cc|}
       \hline
         & & $\rho_1=0.2$  &$\rho_1=0$  &    &  & $\rho_1=0.2$   & $\rho_1=0$  \\

        $p,n$  & $n-p$& $\rho_2=-0.3$  &$\rho_2=0.3$  &  $p,n$  & $n-p$  & $\rho_2=-0.3$   & $\rho_2=0.3$  \\
       \hline
       2,32  & 30 & 0.2657 & 0.1892  & 8,32 & 24 & 0.1202 & 0.1822\\
       \hline
       16,32 & 16 & 0.0519 &0.0281 & 24,32 & 8 & 0.0202 & 0.0198\\
       \hline
       2,64 & 62 & 0.5721  &0.3981 & 32,64 & 32& 0.1190 & 0.1998\\
       \hline
       32,128 &96 & 0.6738 &0.5285 & 64,128& 64 & 0.2853 & 0.1757\\
       \hline
       32,256 & 224 & 0.9291&0.8898 & 64,256 & 192 & 0.9034 & 0.7370\\
       \hline
       128,256 & 128 & 0.6320&0.4225& 16,512 & 496 & 1.0000 & 0.9998\\
       \hline
       64,512 & 448 & 1.0000&0.9989&  128,512 & 384 & 0.9989 & 0.9893\\
       \hline
       256,512 & 256 & 0.9137&0.7530 & 64,1028 & 964 & 1.0000 & 1.0000\\
       \hline
       128,1028 & 900 & 1.0000&1.0000 & 256,1028 &772 & 1.0000 & 1.0000\\
       \hline
     \end{tabular}
  \caption{Empirical power under Gamma(4,1/2) error assumption}\label{table8}
  \end{table}

 \subsection{Parameter estimation under the null hypothesis}

 In practice, if the error terms are not Gaussian, we need to estimate the fourth-order cumulant to perform the test. We now give a suggested estimate under the additional assumption that the error terms are independent under the null hypothesis.
Note that an unbiased estimate of variance $\sigma^2$ under the null hypothesis is
 $$\hat \sigma_n^2=\frac{\gamma_0}{n-p},$$
 and
 \begin{align*}
E\sum_{i=1}^ne_i^4&=3\sigma^4\sum_{i=1}^n \sum_{j_1,j_2}h_{ij_1}^2h_{ij_2}^2+\nu_4\sigma^4\sum_{i=1}^n\sum_{j=1}^nh_{ij}^4=3\sigma^4\tr  (\mb R \circ \mb R)
  + \nu_4\sigma^4\tr(\mb R\du \mb R)^2.
\end{align*}
Then, $\nu_4$ can be estimated by a consistent estimator
 $$
 \hat \nu_4= \frac{\sum_{i=1}^n e_i^4-3\hat\sigma_n^4\tr  (\mb R \circ \mb R)}{\hat\sigma_n^4\tr(\mb R\du \mb R)^2}.
 $$

\section{A general joint CLT for several general quadratic forms}\label{clt}
In this section, we establish a general joint CLT for several general quadratic forms, which helps us to find the asymptotic distribution of the statistics for testing the series correlations. We believe that the result presented below may have its own interest.
\subsection{A brief review of random quadratic forms}
Quadratic forms play an important role not only in mathematical statistics  but also in many other branches of mathematics, such as number theory, differential geometry, linear algebra and differential topology. Suppose ${{\bs\varepsilon}}=\(\varepsilon_1,\varepsilon_2,\cdots,\varepsilon_n\)'$, where $\{\varepsilon_i\}_{i=1}^n$ is a sample of size $n$ drawn from a certain standardized population.
Let $\mb A=(a_{ij})_{n\times n}$ be a matrix. Then, ${{\bs\varepsilon}}'\mb A{{\bs\varepsilon}}=\sum_{i,j}a_{ij}\varepsilon_i\varepsilon_j$ is called a random quadratic form in $\bs\varepsilon$. The random quadratic forms of normal variables, especially when $\mb A$ is symmetric, have been considered by many authors, who have achieved fruitful results. We refer the reader to \cite{Bartlett1960The,Darroch1961Computing,Gart1970Notes,Hsu1999Miscellanea,Forchini2002The,Dik2010The,Al2016On}. Furthermore, many authors have considered the more general situation, where $\bs\varepsilon$ follow a non-Gaussian distribution. For the properties of those types of random quadratic forms, we refer the reader to \cite{Fox1985Central,Cambanis1985Convergence,de1987central,Gregory1995Random,gotze1999asymptotic,liu2009new,Deya2014Invariance,Oliveira2016The} and the references therein.

However, few studies have considered the joint distribution of several quadratic forms. Thus, in this paper, we want to establish a general joint CLT for several random quadratic forms with general distributions.
\subsection{Assumptions and results}
To this end, suppose $$\mb \Xi=
\left(
\begin{array}{cccc}
\varepsilon_{1}^{(1)} & \varepsilon_{2}^{(1)} & \cdots & \varepsilon_{n}^{(1)} \\
\varepsilon_{1}^{(2)} & \varepsilon_{2}^{(2)} & \cdots & \varepsilon_{n}^{(2)} \\
\vdots & \ddots &  & \vdots \\
\varepsilon_{1}^{(q)} & \varepsilon_{2}^{(q)} & \cdots & \varepsilon_{n}^{(q)} \\
\end{array}
\right)
$$
is a $q\times n$ random matrix. Let $\{\mb A_l=a_{i,j}^{(l)}\}, 1\leq i,j\leq n, 1\leq l\leq q$ be $q$ nonrandom $n$-dimensional matrices. Define $Q_l=\sum_{i,j=1}^na_{i,j}^{(l)}\varepsilon_i^{(l)}\varepsilon_j^{(l)}$ for $1\leq l\leq q.$ We are interested in the asymptotic distribution, as $n\to\infty$, of the random vector $\(Q_1,Q_2,\cdots,Q_q\)$, which consists of $q$ random quadratic forms.
Now, we make the following assumptions.
\begin{itemize}
\item[(a)] $\{\varepsilon_j^{(i)}\}_{1\le j\le n,1\leq i\leq q}$ are standard random variables (mean zero and variance one) with uniformly bounded fourth-order moments $M_{4,j}^{(i)}$.
\item[(b)] The columns of $\mb \Xi$ are independent.
\item[(c)] The spectral norms  of the $q$ $n \times n$ square matrices $\{\mb A_l\}_{1\leq l\leq q}$ are uniformly  bounded in $n$.
\end{itemize}
Clearly, for $1\leq l\leq q$, we have $\E Q_l=\tr \mb A_l$, and for $1\leq l_1,l_2\leq q$, we obtain \begin{align}
  &{{\rm Cov}\(Q_{l_1},Q_{l_2}\)}=\E\(Q_{l_1}-\E Q_{l_1}\)\(Q_{l_2}-\E Q_{l_2}\)\\\notag
  =&\E\(\sum_{i=1}^na_{i,i}^{(l_1)}\(\(\varepsilon_i^{(l_1)}\)^2-1\)+\sum_{i\neq j}a_{i,j}^{(l_1)}\varepsilon_i^{(l_1)}\varepsilon_j^{(l_1)}\)\(\sum_{i=1}^na_{i,i}^{(l_2)}\(\(\varepsilon_i^{(l_2)}\)^2-1\)+\sum_{i\neq j}a_{i,j}^{(l_2)}\varepsilon_i^{(l_2)}\varepsilon_j^{(l_2)}\)\\\notag
  =&\sum_{i=1}^na_{i,i}^{(l_1)}a_{i,i}^{(l_2)}\E\(\(\varepsilon_i^{(l_1)}\)^2-1\)\(\(\varepsilon_i^{(l_2)}\)^2-1\)
  +\sum_{i\neq j}^na_{i,j}^{(l_1)}a_{i,j}^{(l_2)}\E\varepsilon_i^{(l_1)}\varepsilon_i^{(l_2)}\E\varepsilon_j^{(l_1)}\varepsilon_j^{(l_2)}\\\notag
  &+\sum_{i\neq j}^na_{i,j}^{(l_1)}a_{j,i}^{(l_2)}\E\varepsilon_i^{(l_1)}\varepsilon_i^{(l_2)}\E\varepsilon_j^{(l_1)}\varepsilon_j^{(l_2)}.
\end{align}
Let $l_1=l_2=l$; then, we have $${\rm Var}(Q_l)=\sum_{i=1}^n(M_{4,i}^{(l)}-3)\(a_{i,i}^{(l)}\)^2+\tr \mb A_l\mb A_l'+\tr \mb A_l^2.$$
Thus, according to assumptions $(a)-(c)$, for any $1\leq l\leq q,$ ${\rm Var}(Q_l)$ at most has the same order as $n$. This result also holds for any ${{\rm Cov}\(Q_{l_1},Q_{l_2}\)}$ by applying the Cauchy-Schwartz inequality.
We then have the following theorem.
\begin{thm}\label{th:CLT}
 In addition to assumptions (a)-(c), suppose that there exists an $i$ such that ${\rm Var}(Q_i)$ has the same order as $n$ when $n\to \infty$.
  Then, the distribution of the random vector $n^{-1/2}(Q_1, Q_2,\cdots,Q_q)$
  is asymptotically $q$-dimensional normal.
\end{thm}

\subsection{Proof of Theorem~\protect\ref{th:CLT}}

We are now in position to present the proof of the joint CLT via the method of moments. The procedure of the proof is similar to that in \cite{yin2018test} but is more complex since we need to establish the CLT for a $q$-dimensional, rather than 2-dimensional, random vector. Moreover, we do not assume the underlying distribution to be symmetric and identically distributed. The proof is separated into three steps.
\subsubsection{Step 1: Truncation}

Noting that $\sup_{i,j}\E(\varepsilon_{j}^{(i)})^4<\infty$, $j=1,\cdots,n, i=1,\cdots,q$, for any $\delta>0$, we have
$\sup_{i,j} \delta^{-4}nP(|(\varepsilon_{j}^{(i)})|\ge \delta n^{1/4})\to 0.$ Thus, we may select a sequence $\delta_n\to 0$ such that
$\sup_{i,j}\delta_n^{-4}nP(|(\varepsilon_{j}^{(i)})|>\delta_n n^{1/4})\to 0$. The convergence rate of $\delta_n $ to 0 can be made arbitrarily slow.  Define $(\widetilde Q_1,\widetilde Q_2,\cdots,\widetilde Q_q)$ to be the analogue of $(Q_1,Q_2,\cdots,Q_q)$ with $\varepsilon_{j}^{(i)}$
replaced by $\tilde\varepsilon_{j}^{(i)}$,
where $\tilde \varepsilon_{j}^{(i)}=\varepsilon_{j}^{(i)}I(|(\varepsilon_{j}^{(i)})|<\delta_n n^{1/4})$. Then,
$$
P\((\widetilde Q_1,\widetilde Q_2,\cdots,\widetilde Q_q)\ne (Q_1,Q_2,\cdots,Q_q)\)\le \sum_{i=1}^q\sum_{j=1}^nP(\varepsilon_{j}^{(i)}\ne \tilde\varepsilon_{j}^{(i)})\leq qnP(|(\varepsilon_{j}^{(i)})|\ge \delta_n n^{1/4})\to 0.
$$
Therefore, we need only to investigate the limiting distribution of the vector $(\widetilde Q_1,\widetilde Q_2,\cdots,\widetilde Q_q)$.

\subsubsection{Step 2: Centralization and Rescaling}

Define $(\breve Q_1,\breve Q_2,\cdots,\breve Q_q)$ to be the analogue of $(\widetilde Q_1,\widetilde Q_2,\cdots,\widetilde Q_q)$ with $\widetilde \varepsilon_{j}^{(i)}$
replaced by $\breve\varepsilon_{j}^{(i)}=\frac{\tilde\varepsilon_{j}^{(i)}-\E \tilde\varepsilon_{j}^{(i)}}{\sqrt{{\rm Var}\tilde\varepsilon_{j}^{(i)}}}$.
Denote by $d(X,Y)=\sqrt{E|X-Y|^2}$ the distance between two  random variables $X$ and $Y$. Additionally, denote $\widetilde\varepsilon^{(i)}=\(\widetilde\varepsilon_{1}^{(i)},\widetilde\varepsilon_{2}^{(i)},\cdots,\widetilde\varepsilon_{n}^{(i)}\)'$, $\breve\varepsilon^{(i)}=\(\breve\varepsilon_{1}^{(i)},\breve\varepsilon_{2}^{(i)},\cdots,\breve\varepsilon_{n}^{(i)}\)'$ and $\Psi^{(i)}={\rm Diag}\({\sqrt{{\rm Var}\tilde\varepsilon_{1}^{(i)}}},{\sqrt{{\rm Var}\tilde\varepsilon_{2}^{(i)}}},\cdots,{\sqrt{{\rm Var}\tilde\varepsilon_{n}^{(i)}}}\).$
We obtain that for any $l=1,\cdots,q,$
\begin{align}\label{eqbb1}
d^2(\widetilde Q_l,\breve Q_l)&={\E|\( {\breve \varepsilon}'^{(l)}\mb A_l \breve\varepsilon^{(l)}-\widetilde\varepsilon'^{(l)}\mb A_l \widetilde\varepsilon^{(l)}\)|^2}\\\notag
&={\E|\( {\breve\varepsilon}'^{(l)}\mb A_l \breve\varepsilon^{(l)}-\breve\varepsilon'^{(l)}\Psi^{(l)}\mb A_l \Psi^{(l)}\breve\varepsilon^{(l)}-{\E \tilde\varepsilon}'^{(l)}\Psi^{(l)}\mb A_l \Psi^{(l)}\E \tilde\varepsilon^{(l)}\)|^2}\\\notag
&\leq 2\(\E| {\breve\varepsilon}'^{(l)}\mb A_l \breve\varepsilon^{(l)}-\breve\varepsilon'^{(l)}\Psi^{(l)}\mb A_l \Psi^{(l)}\breve\varepsilon^{(l)}|^2+|{\E \tilde\varepsilon}'^{(l)}\Psi^{(l)}\mb A_l \Psi^{(l)}\E \tilde\varepsilon^{(l)}|^2\)\\\notag
&\triangleq2\(\Upsilon_{l,1}+\Upsilon_{l,2}\).\\\notag
\end{align}
Noting that $\breve\varepsilon_{j}^{(i)}$'s are independent random variables with 0 means and unit variances, it follows that
\begin{align*}
\Upsilon_{l,1}&=\E\( {\breve\varepsilon}'^{(l)}\mb A_l \breve\varepsilon^{(l)}-\breve\varepsilon'^{(l)}\Psi^{(l)}\mb A_l \Psi^{(l)}\breve\varepsilon^{(l)}\)^2=\E\( {\breve\varepsilon}'^{(l)}\(\mb A_l-\Psi^{(l)}\mb A_l \Psi^{(l)}\)\breve\varepsilon^{(l)}\)^2\\\notag
&=\nu_4\tr\(\(\mb A_l-\Psi^{(l)}\mb A_l \Psi^{(l)}\)\du\(\mb A_l-\Psi^{(l)}\mb A_l \Psi^{(l)}\)\)+\tr\(\mb A_l-\Psi^{(l)}\mb A_l \Psi^{(l)}\)^2\\\notag
&+\tr\(\(\mb A_l-\Psi^{(l)}\mb A_l \Psi^{(l)}\)\(\mb A_l-\Psi^{(l)}\mb A_l \Psi^{(l)}\)'\).
\end{align*}
Since $\E\tilde\varepsilon_j^{(l)}=\E\varepsilon_j^{(l)}I\(|\varepsilon_j^{(l)}|\geq\delta_nn^{1/4}\)\leq C\delta_n^{-3}n^{-3/4},$ and $$1-\E(\tilde\varepsilon_j^{(l)})^2=\E(\varepsilon_j^{(l)})^2I\(|\varepsilon_j^{(l)}|\geq\delta_nn^{1/4}\)\leq C\delta_n^{-2}n^{-1/2},$$
we know that
\begin{align}
\|\(\mb I-\Psi^{(l)}\)\|=\max_{j=1,\cdots,n}{|1-\sqrt{{\rm Var}\tilde\varepsilon_j^{(l)}}}|\leq \max_{j=1,\cdots,n}\sqrt{|1-{\rm Var}\tilde\varepsilon_j^{(l)}}|\leq C\delta_n^{-1}n^{-1/4}.
\end{align}
Then, we have $$\|\(\mb A_l-\Psi^{(l)}\mb A_l \Psi^{(l)}\)\|\leq \|\mb A_l\|\|\(\mb I-\Psi^{(l)}\)\|+\|\(\mb I-\Psi^{(l)}\)\|\|\mb A_l\|\|\Psi^{(l)}\|=O(\delta_n^{-1}n^{-1/4}).$$
It follows that $\Upsilon_{l,1}=O(\delta_n^{-2}n^{1/2})$ and $\Upsilon_{l,2}\leq\|\Psi^{(l)}\|^2\|\mb A_l\|\sum_{j=1}^n\(\E\tilde\varepsilon_j^{(l)}\)^2=O(\delta_n^{-6}n^{-1/2}).$
By combining the above estimates, we obtain that  $d(\widetilde Q_l,\breve Q_l)=O(\delta_n^{-1}n^{1/4})$ for $l=1,\cdots,q.$

Noting that the entries in the covariance matrix of the random vector $(Q_1,Q_2,\cdots,Q_q)'$ have at most the same order as $n$, we conclude that $n^{-1/2}(Q_1,Q_2,\cdots,Q_q)'$ has the same limiting distribution as the random vector $n^{-1/2}(\breve Q_1,\breve Q_2,\cdots,\breve Q_q)'$. Therefore, we shall subsequently assume that  $|\varepsilon_{j}^{(i)}|\le \delta_n n^{1/4}$ holds in the proof of the CLT.

\subsubsection{Step 3: Completion of the proof}

Let $\alpha_1, \cdots, \alpha_q$ be q real numbers satisfying $\sum_{l=1}^q\alpha_l^2\neq 0$.
We show that for any $k$,
\begin{align}\label{eqclt1}
  &\E\(\sum_{l=1}^q\alpha_l \(Q_l-\E Q_l\)\)^k\\\notag
  =&\begin{cases}(k-1)!!\(\sum_{l_1=1}^q\sum_{l_2=1}^q\alpha_{l_1}\alpha_{l_2}{\rm Cov}\(Q_{l_1},Q_{l_2}\)\)^{k/2}(1+o(1))& \mbox{for $k$ is even,}\\
    o(n^{k/2})&\mbox{for $k$ is odd}.\end{cases}
\end{align}

\begin{figure}[H]
  \includegraphics[width=0.8\textwidth]{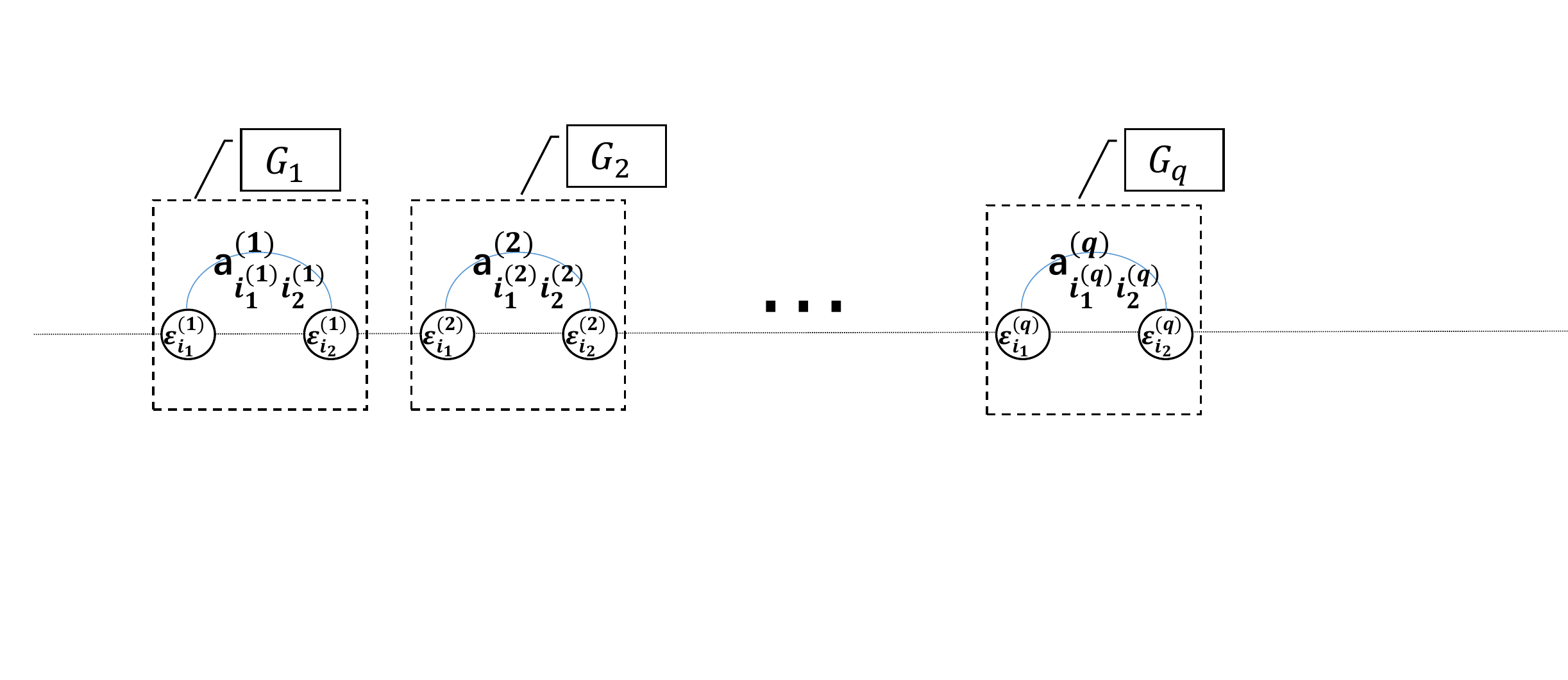}
  \caption{The basic graphs.}\label{basicg}	
\end{figure}

Write
\begin{align}
  \E\(\sum_{l=1}^q\alpha_l \(Q_l-\E Q_l\)\)^k
  =\sum_{\substack{k_1,k_2,\dots,k_q\geq 0\\\sum_{l=1}^q{k_{l}}=k}}\frac{k!}{\prod_{l=1}^qk_l!}\prod_{l=1}^q\alpha_l^{k_l}
  \E\prod_{l=1}^q(Q_l-\E Q_l)^{k_l}.\label{eqclt2}
\end{align}
Draw a parallel line and for given $2q$ numbers $i_1^{(1)},i_2^{(1)},\cdots,i_{1}^{(q)},i_{2}^{(q)}$ on this line, draw q simple graphs $$G_1(i_1^{(1)},i_2^{(1)}), G_2(i_1^{(2)},i_2^{(2)}),\cdots,G_q(i_1^{(q)},i_2^{(q)})$$ from $i_1^{(l)}$ to $i_2^{(l)}$ for $l=1,\cdots,q$. For any $l=1,\cdots,q,$ we use the edge $(i_1^{(l)},i_2^{(l)})$ to indicate that the entry lies in the $i_1^{(l)}$-th row and $i_2^{(l)}$-th column of matrix $\mb A_l$, denoted as $a_{i_1^{(l)},i_2^{(l)}}^{(l)}$. The two vertices $i_1^{(l)}$ and $i_2^{(l)}$ correspond to random variables $\varepsilon_{i_1}^{(l)}$ and $\varepsilon_{i_2}^{(l)}$, respectively.
Thus, the graph $G_l(i_1^{(l)},i_2^{(l)})$ corresponds to $A_{G_l(i_1^{(l)},i_2^{(l)})}\varepsilon_{G_l(i_1^{(l)},i_2^{(l)})}$, where $A_{G_l(i_1^{(l)},i_2^{(l)})}=a_{i_1^{(l)},i_2^{(l)}}^{(l)}$ and $\varepsilon_{G_l(i_1^{(l)},i_2^{(l)})}=\varepsilon_{i_1}^{(l)}\varepsilon_{i_2}^{(l)}$.
We call $G_1,\cdots,G_q$ the basic graphs. Figure\ref{basicg} shows the basic graphs.

\begin{figure}[H]
  \includegraphics[width=0.8\textwidth]{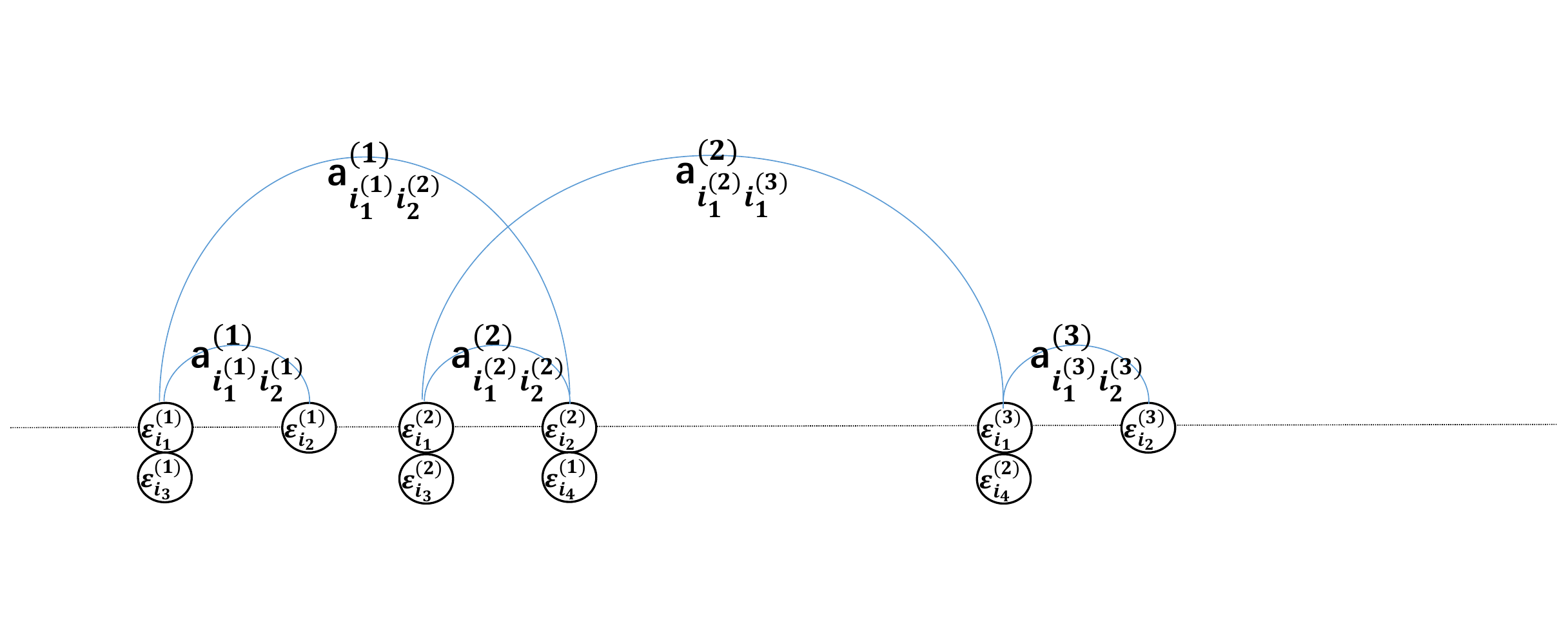}
  \caption{A graph associated with terms that satisfy $i_1^{(1)}$ equal to $i_3^{(1)}$, $i_1^{(2)}$ equal to $i_3^{(2)}$, $i_2^{(2)}$ equal to $i_4^{(1)}$ and $i_1^{(3)}$ equal to $i_4^{(2)}$.}\label{terms}	
\end{figure}

Now, we draw $k_l$ $G_l$ graphs for $l=1,\cdots,q$ and denote them by $G_{l,\ell}$,$\ell=1,\cdots,k_l$.
Note that
\begin{align}
  &\E\prod_{l=1}^q(Q_l-\E Q_l)^{k_l}=\E\prod_{l=1}^q\(\sum_{i_{1}^{(l)},i_{2}^{(l)}}a_{i_{1}^{(l)},i_{2}^{(l)}}^{(l)}\(\varepsilon_{i_{1}^{(l)}}\varepsilon_{i_{2}^{(l)}}-\E\varepsilon_{{i_{1}^{(l)}}}\varepsilon_{i_{2}^{(l)}}\)\)^{k_l}\\\notag
  =&\sum_{G_{l,\ell}}\E\(
    \prod_{l=1}^q\prod_{\ell=1}^{k_l}(A_{G_{l,\ell}}\varepsilon_{{G_l,\ell}}-\E A_{G_{l,\ell}}\varepsilon_{{G_{l,\ell}}})\),
\end{align}
where the summation runs over all possibilities of the $G_1,\cdots,G_q$ graphs (according to the values of $k_1,\cdots,k_l$ and $i_1^{(1)},\cdots,i_{2,k_q}^{(q)}$). We have now completed the step of associating the terms in the expression of $\E\prod_{l=1}^q(Q_l-\E Q_l)^{k_l}$ with graphs. For example, for $q=4,k=5,k_1=2.k_2=2,k_3=1,k_4=0$, the graph in Figure \ref{terms} is associated with the term
\begin{align}
  \sum_{i_1^{(1)},i_2^{(1)},i_1^{(2)},i_2^{(2)},i_1^{(3)},i_2^{(3)}=1}^n&a_{i_1^{(1)}i_2^{(1)}}^{(1)}a_{i_1^{(1)}i_2^{(2)}}^{(1)}a_{i_1^{(2)}i_2^{(2)}}^{(2)}
  a_{i_1^{(2)}i_1^{(3)}}^{(2)}a_{i_1^{(3)}i_2^{(3)}}^{(3)}\(\varepsilon_{i_1^{(1)}}\varepsilon_{i_2^{(1)}}-\E\varepsilon_{i_1^{(1)}}\varepsilon_{i_2^{(1)}}\)\\\notag
  &\times
  \(\varepsilon_{i_1^{(1)}}\varepsilon_{i_2^{(2)}}-\E\varepsilon_{i_1^{(1)}}\varepsilon_{i_2^{(2)}}\)
  \(\varepsilon_{i_1^{(2)}}\varepsilon_{i_2^{(2)}}-\E\varepsilon_{i_1^{(2)}}\varepsilon_{i_2^{(2)}}\)\\\notag
  &\times
  \(\varepsilon_{i_1^{(2)}}\varepsilon_{i_1^{(3)}}-\E\varepsilon_{i_1^{(2)}}\varepsilon_{i_1^{(3)}}\)
  \(\varepsilon_{i_1^{(3)}}\varepsilon_{i_2^{(3)}}-\E\varepsilon_{i_1^{(3)}}\varepsilon_{i_2^{(3)}}\)
  .
\end{align}

We next classify all the terms in the above summation into three groups. Group one contains all the terms whose corresponding combined graph $G$ has at least one subgraph that does not have any vertices  coincident  with vertices of the other subgraphs. Group two contains all the terms whose corresponding combined graph $G$ has at least one vertex that is not coincident with any other vertices. All the other terms are classified into the third group. Since the entries in $\varepsilon$ are independent with 0 means, all the terms in group one and group two are equal to 0.

\begin{figure}[H]
  \includegraphics[width=0.8\textwidth]{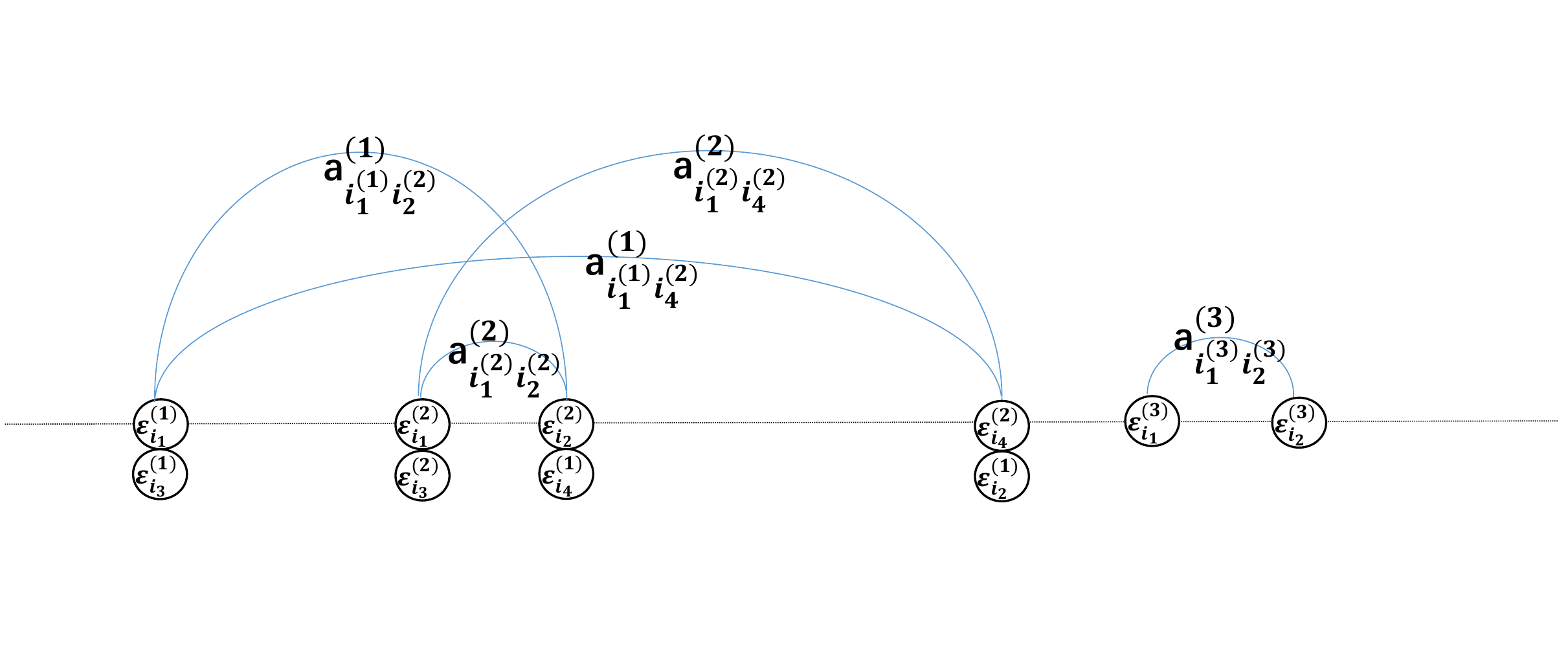}
  \caption{A graph associated with terms classified into group one.}\label{terms2}	
\end{figure}

For example, for $q=4,k=5,k_1=2.k_2=2,k_3=1,k_4=0$, the term associated with the graph shown in Figure \ref{terms} is classified into group two, while the term associated with the graph shown in Figure \ref{terms2} is classified into group one.

\begin{figure}[H]
  \includegraphics[width=0.8\textwidth]{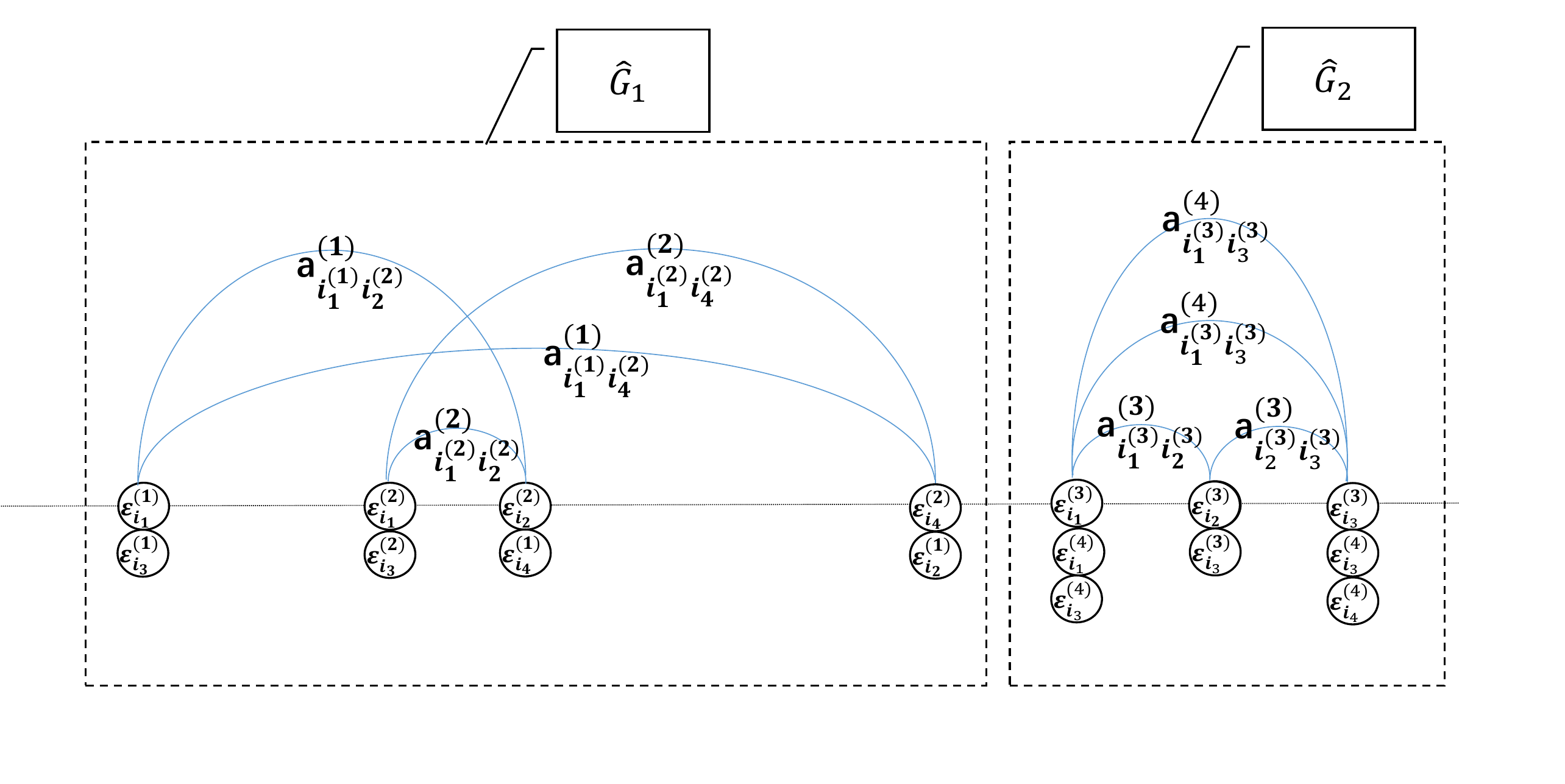}
  \caption{A graph G associated with terms classified into group three.}\label{terms3}	
\end{figure}

Therefore, we need only to evaluate the sum of terms that belong to the third group.
Suppose the combined graph  $G$ contains
$\pi$ connected pieces $\hat G_1,\cdots,\hat G_\pi$ consisting of
$\phi_1,\cdots,\phi_{\pi}$ subgraphs ($G_1,\cdots,G_q$). Clearly,
$\phi_1,\phi_2,\cdots,\phi_{\pi}\ge 2$ since any subgraph must have at least one  vertex  coincident  with vertices of the other subgraphs; hence, $\pi\le k/2$ since $\phi_1+\cdots+\phi_\pi=k$. For example, in graph G shown in Figure \ref{terms3}, $q=4$, $k=8$, $\pi=2$, $\phi_1=\phi_2=4$.

We then introduce some necessary definitions and lemmas about graph-associated multiple matrices for the purpose of calculating the contributions of those terms in group three.

We first give two definitions:
\begin{definition}[two-edge connected]
  A graph $G$ is called {\bf two-edge connected} if the resulting subgraph is still connected after removing any edge from G.
\end{definition}

\begin{definition}[cutting edge]
  An edge $e$ in a graph $G$ is called a {\bf cutting edge} if deleting this edge results in a disconnected subgraph.
\end{definition}

\begin{figure}[H]
  \includegraphics[width=0.8\textwidth]{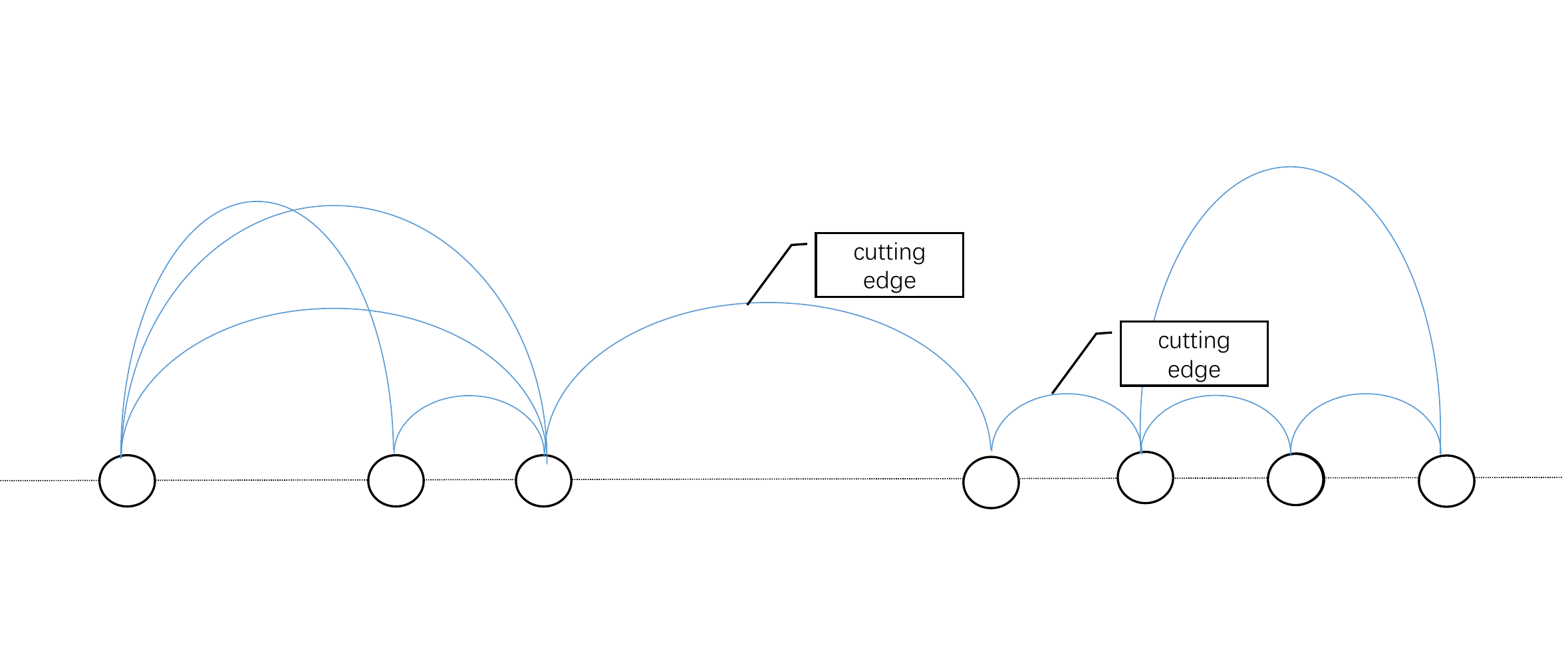}
  \caption{An example of graph that is not two-edge connected. There are two cutting edges.}\label{two-edge2}	
\end{figure}

Clearly, a graph is a two-edge connected graph if and only if there is no cutting edge. \ref{two-edge1} below shows an example of a two-edge connected graph, while the graph show in Figure \ref{two-edge2} is not a two-edge connected graph and has two cutting edges.

\begin{figure}[H]
  \includegraphics[width=0.8\textwidth]{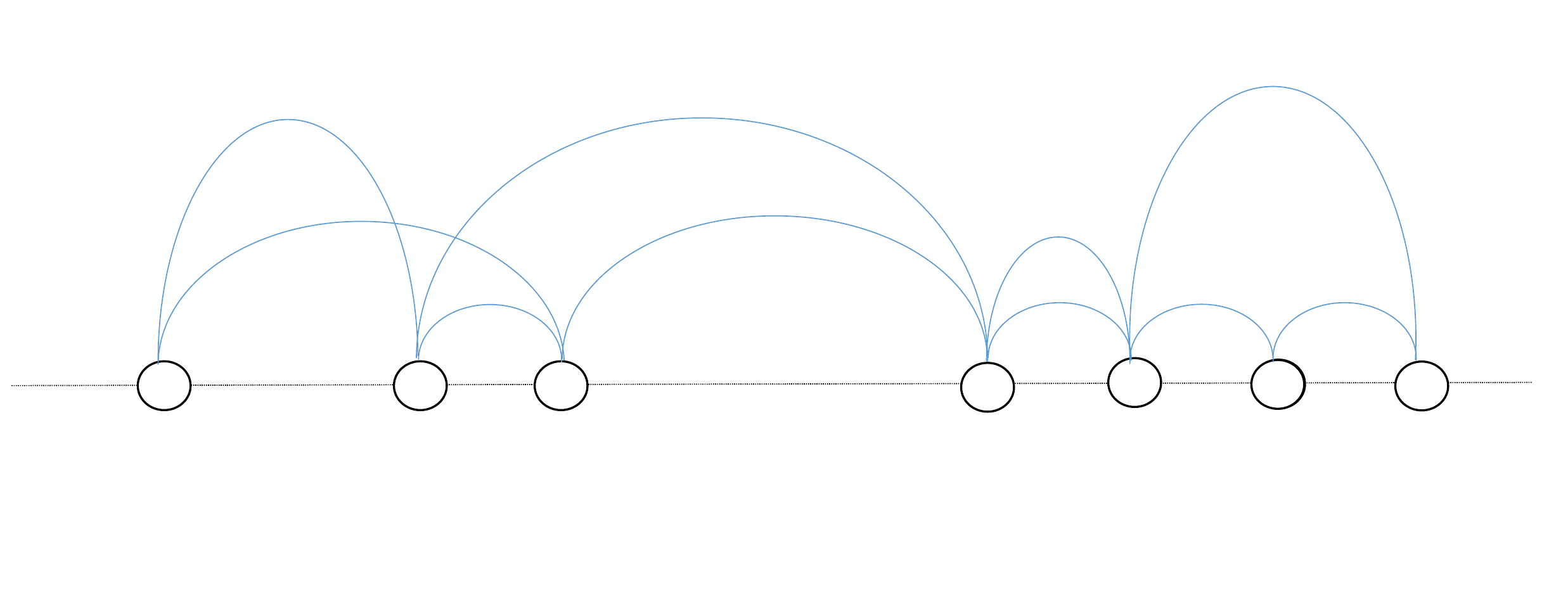}
  \caption{An example of a two-edge connected graph.}	\label{two-edge1}	
\end{figure}

Now, we shall introduce the following lemma.
\begin{lemma}\label{lm2}
 Suppose that $\mb G=\(\mb V,\mb E, \mb F\)$ is a two-edge connected graph with $t$ vertices and $k$ edges. Each vertex $i$ corresponds to an integer $m_i \geq 2$, and each edge $e_j$ corresponds to a matrix $\mb T^{(j)}=\(t_{\alpha,\beta}^{(j)}\),\ j=1,\cdots,k$ with consistent dimensions, that is, if $F(e_j)=(f_i(e_j),f_e(e_j))=(g,h),$ then the matrix $\mb T^{\(j\)}$ has dimensions $m_g\times m_h$. Define $\mb v=(v_1,v_2,\cdots,v_t)$ and
 \begin{align}
   T'=\sum_{\mb v}\prod_{j=1}^kt_{v_{f_i(e_j)},v_{f_e(e_j)}}^{(j)},
 \end{align}
 where the summation $\sum_{\mb v}$ is taken for $v_i=1,2,\cdots, m_i, \ i=1,2,\cdots,t.$ Then, for any $i\leq t$,
 we have
 $$|T'|\leq m_i\prod_{j=1}^k\|\mb T^{(j)}\|.$$
\end{lemma}
For the proof of this lemma, we refer the reader to section $\mb A.4.2$ in \cite{bai2010spectral}.

Now, suppose that the connected piece $\hat G_{\varphi}$ ($1\leq \varphi \leq \pi$) consists of $\phi_{\varphi}$ subgraphs ($G_1,\cdots,G_q$). Then, the number of edges in $\hat G_{\varphi}$ is exactly $\phi_{\varphi}$. Let $\upsilon_{\varphi}$ denote the number of noncoincident vertices (in graph G shown in Figure \ref{terms3}, $\upsilon_{1}=4, \upsilon_{2}=3$). Denote those vertices by $V_{\phi,1},\cdots,V_{\phi,\upsilon_{\varphi}}$. Additionally, denote the degree of those vertices by $\omega_{\phi,1},\cdots,\omega_{\phi,\upsilon_{\varphi}}$. Clearly, $\upsilon_{\varphi}\leq \phi_{\varphi}$ since the total degree is $2\phi_{\varphi}$ and the degrees of all vertices are at least 2.

Note that $\E\prod_{t=1}^{\phi_{\varphi}}(\mb A_{G_{\varphi t}}\bs\varepsilon_{G_{\varphi t}}-\E\mb A_{G_{\varphi t}}\bs\varepsilon_{G_{\varphi t}})
=\mb A_{G_\varphi}\E\prod_{t=1}^{\phi_{\varphi}}(\bs\varepsilon_{G_{\varphi t}}-\E\bs\varepsilon_{G_{\varphi t}})$.
We now focus on estimating the relationship between $\upsilon_{\varphi}$ and $\phi_{\varphi}$.
\begin{itemize}
  \item Case (1): If $\upsilon_{\varphi}=\phi_{\varphi}$, then all the vertices in $\hat G_{\varphi}$ are of degree 2; thus, $\hat G_{\varphi}$ is an Euler graph, which is a circle and is therefore two-edge connected. It follows from Lemma \ref{lm2} that $\sum_{\hat G_{\varphi}}\mb A_{\hat G_{\varphi}}=O(n).$ Since the fourth moment of the underlying distribution is finite, we have $|E\varepsilon_{\hat G_{\varphi}}|=O(1)$.
      An example of a graph in this case with $\phi_{\varphi}=8$ is shown in Figure \ref{case1}.

\begin{figure}[H]
  \includegraphics[width=0.8\textwidth]{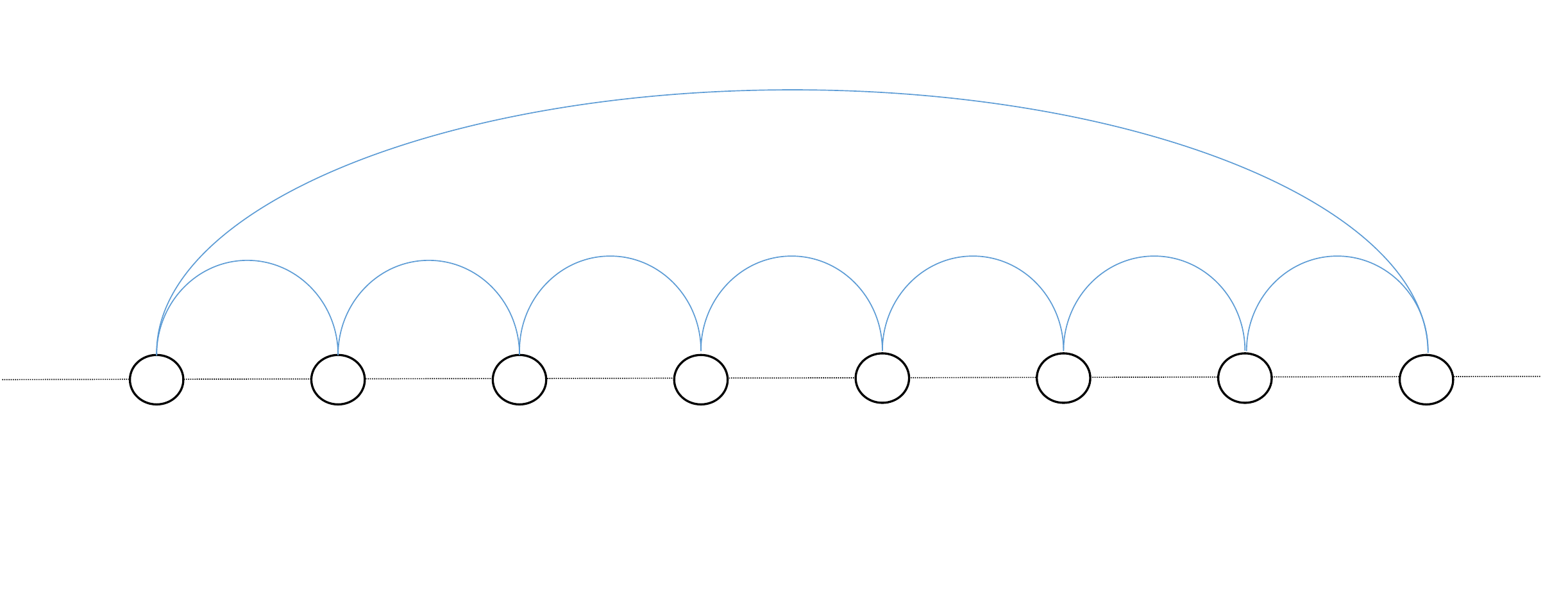}
  \caption{An example of graph that falls into Case (1). }\label{case1}	
\end{figure}

  \item Case (2): If there is exactly two vertices of degree 3 and all other vertices are of degree 2, then the two vertices of degree 3 must lie on the two ``sides" of $\hat G_{\varphi}.$ There are two types of graphs that satisfy these conditions, as shown in Figure \ref{case2}.
\begin{figure}[H]
  \includegraphics[width=0.8\textwidth]{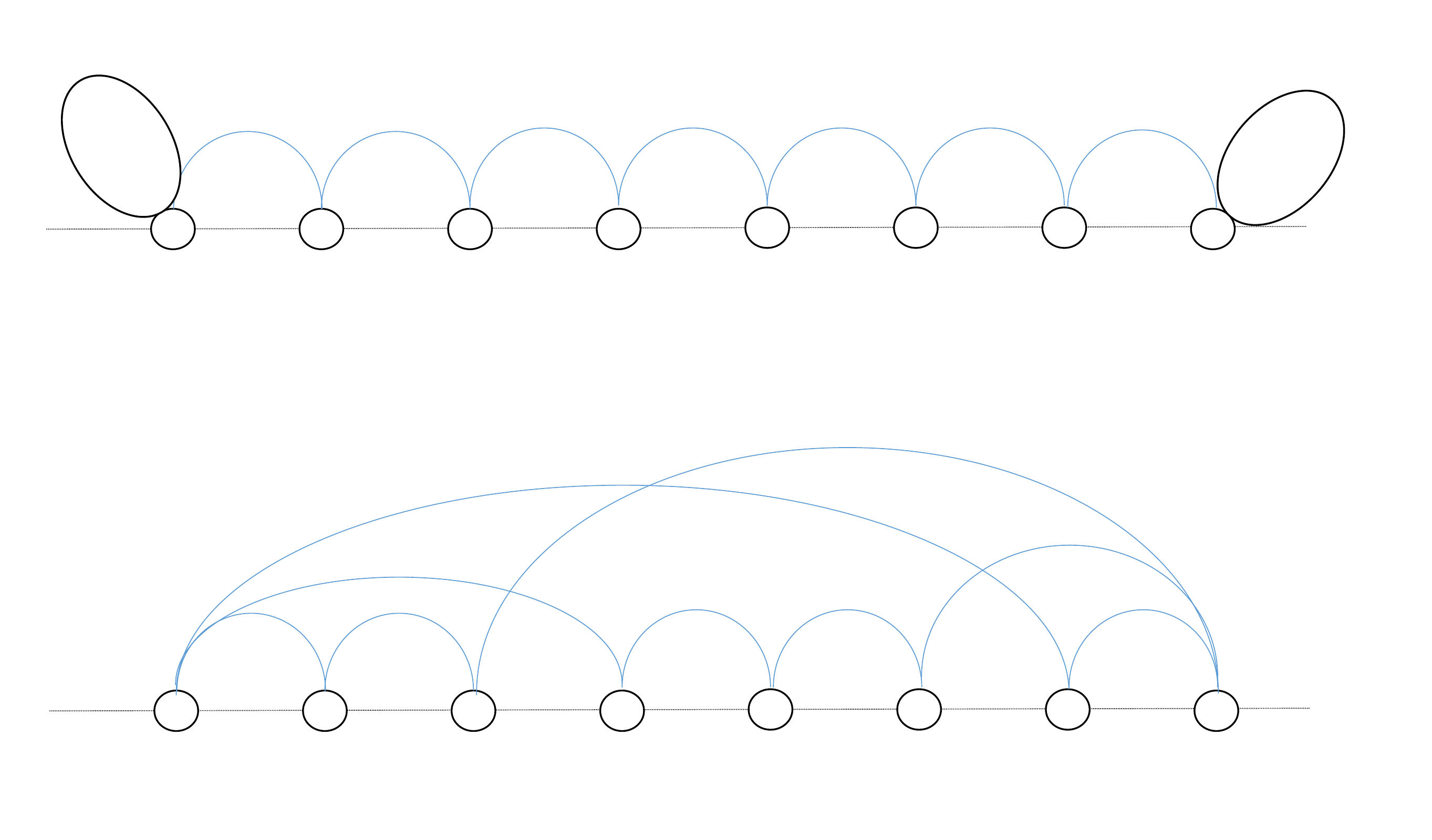}
  \caption{Two types of graphs that fall into Case (2). }\label{case2}	
\end{figure}
      All graphs of the second type are clearly two-edge connected. For the first type of graph, we have $\sum_{\hat G_{\varphi}}\mb A_{\hat G_{\varphi}}=\sum_{i,j=1}^n b_{i,j},$ where $\mathfrak{B}=(b_{i,j})_{n\times n}=\mathfrak{D_1}\mathfrak{A}\mathfrak{D_2}$ with $\mathfrak{D_1}$ and $\mathfrak{D_2}$ being diagonal matrices with a bounded spectrum norm.
        The above arguments imply that we also have $\sum_{\hat G_{\varphi}}\mb A_{\hat G_{\varphi}}=O(n)$ and $|E\varepsilon_{\hat G_{\varphi}}|=O(1)$.
   \item Case (3): If there is exactly one vertex of degree 4 and all other vertices are of degree 2, then, similarly to Case (1), we have $\sum_{\hat G_{\varphi}}\mb A_{\hat G_{\varphi}}=O(n).$ Moreover, we still have $|E\varepsilon_{\hat G_{\varphi}}|=O(1).$
  \item Case (4): $\hat G_{\varphi}$ is a graph that does not fall into the above three cases. Then, suppose there are $\kappa_{\varphi}$ vertices in  $\hat G_{\varphi}$ with degrees larger than 4. Without loss of generality, denote these vertices as $\varrho_{\phi,1},\cdots,\varrho_{\phi,\kappa_{\varphi}}$. Choose a minimal spanning tree $\hat G_{\varphi}^{0}$ from $\hat G_{\varphi}$. Denote the remaining graph by $\hat G_{\varphi}^{1}$. An example graph $\hat G_{\varphi}$ that falls into this case is shown in Figure \ref{case4}.
       \begin{figure}[H]
  \includegraphics[width=0.8\textwidth]{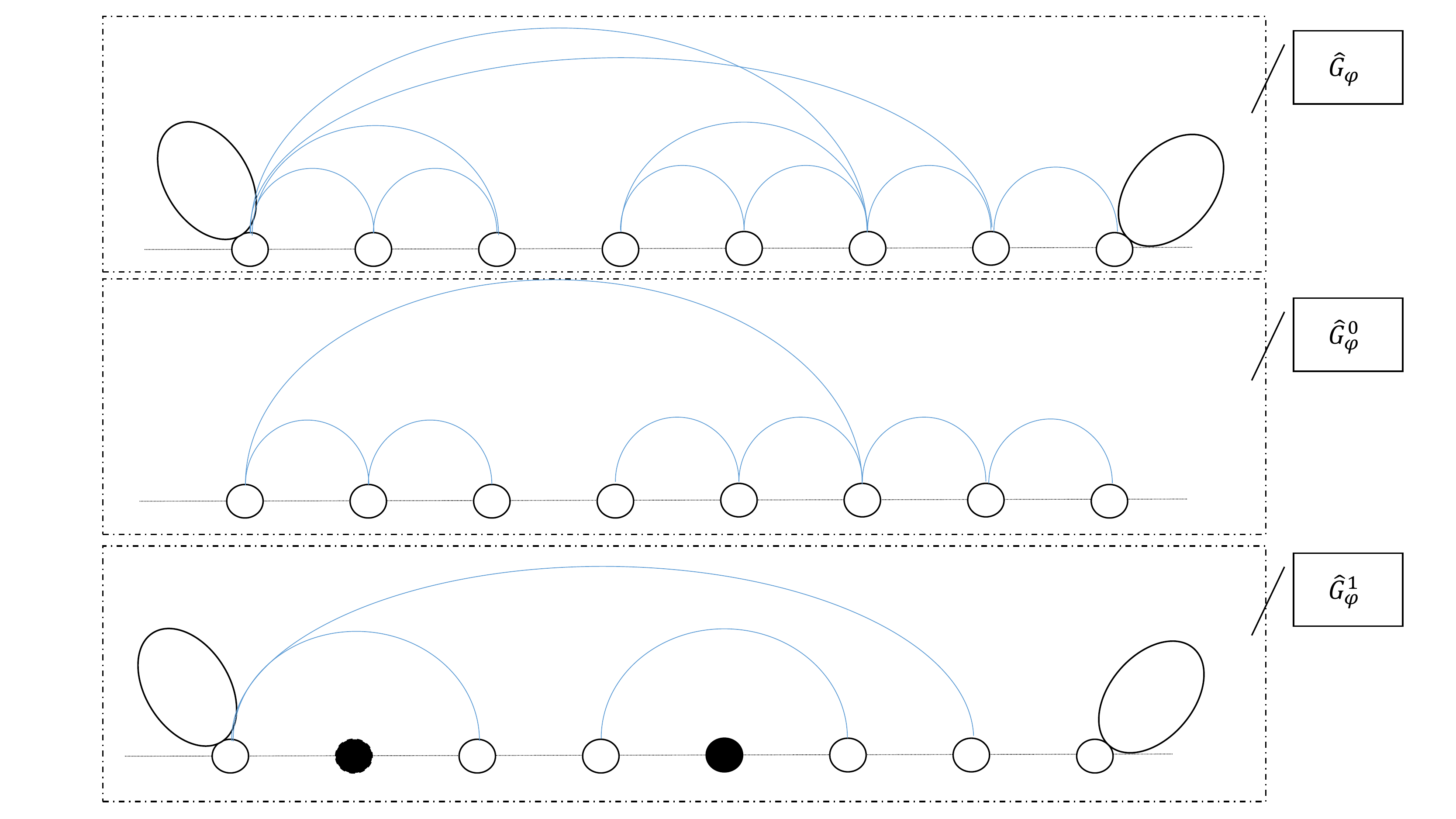}
  \caption{An example graph $\hat G_{\varphi}$ that falls into Case (4). $\hat G_{\varphi}^{0}$ is a minimal spanning tree of $\hat G_{\varphi}$, and $\hat G_{\varphi} ^1$ is the remaining graph.}\label{case4}	
\end{figure}
      Then, by the Cauchy-Schwarz inequality, we have
      $$\sum_{\hat  G_{\varphi}}\mb A_{\hat G_{\varphi}}\leq \(\sum_{\hat G_{\varphi}}\mb A_{\hat G_{\varphi}^{0}\cup \hat G_{\varphi}^{0}}\)^{1/2}\(\sum_{\hat G_{\varphi}}\mb A_{\hat G_{\varphi}^{1}\cup \hat G_{\varphi}^{1}}\)^{1/2}.$$
      Note that all the degrees of the vertices in ${\hat G_{\varphi}^{0}\cup \hat G_{\varphi}^{0}}$ are even; thus, ${\hat G_{\varphi}^{0}\cup \hat G_{\varphi}^{0}}$ is an Euler graph. Additionally, note that $\hat G_{\varphi}^{0}$ contains all the vertices of $\hat G_{\varphi}$. It follows from Lemma \ref{lm2} that
      $\sum_{\hat G_{\varphi}}\mb A_{\hat G_{\varphi}^{0}\cup \hat G_{\varphi}^{0}}=O(n).$ For the same reason, since all the degrees of the vertices in ${\hat G_{\varphi}^{1}\cup \hat G_{\varphi}^{1}}$ are even and the number of disconnected subgraphs (including isolated vertices, if they exist) of ${\hat G_{\varphi}^{1}\cup \hat G_{\varphi}^{1}}$ is at most $\upsilon_{\varphi}$. Thus, we have $\sum_{\hat G_{\varphi}}\mb A_{\hat G_{\varphi}^{0}\cup \hat G_{\varphi}^{0}}=O(n^{\upsilon_{\varphi}}).$ Now, we estimate $|E\varepsilon_{\hat G_{\varphi}}|.$ If $\kappa_{\varphi}=0,$ then we have $|E\varepsilon_{\hat G_{\varphi}}|=O(1).$ If $\kappa_{\varphi} \geq 1,$ denote the degrees of $\varrho_{\phi,1},\cdots,\varrho_{\phi,\kappa_{\varphi}}$ as $d_{\varrho_{\phi,1}},\cdots,d_{\varrho_{\phi,\kappa_{\varphi}}}$, respectively. Since the underlying variables are truncated and $\delta_n\to 0$, we obtain that, for large $n$,
\begin{align}
  |E\varepsilon_{\hat G_{\varphi}}|\le \(\delta_n n^{1/4}\)^{\(\sum_{i=1}^{\kappa_{\varphi}}\(d_{\varrho_{\phi,i}}-4\)\)}=o\( n^{\frac{\sum_{i=1}^{\kappa_{\varphi}}\(d_{\varrho_{\phi,i}}-4\)}{4}}\),
\end{align}
which implies that
\begin{align}
  &\mb A_{\hat G_\varphi}\E\prod_{t=1}^{\phi_{\varphi}}(\bs\varepsilon_{\hat G_{\varphi t}}-\E\bs\varepsilon_{\hat G_{\varphi t}})=o(n^{1/2}n^{{\upsilon_{\varphi}}/2}n^{\frac{\sum_{i=1}^{\kappa_{\varphi}}\(d_{\varrho_{\phi,i}}-4\)}{4}})=o(n^{\frac{{\sum_{i=1}^{\kappa_{\varphi}}d_{\varrho_{\phi,i}}}-4\kappa_{\varphi}+2\upsilon_{\varphi}+2}{4}})\\\notag
  =&o\(n^{{\frac{{\sum_{i=1}^{\kappa_{\varphi}}d_{\varrho_{\phi,i}}}+2\(\upsilon_{\varphi}-\kappa_{\varphi}\)-2\kappa_{\varphi}+2}{4}}}\)=o(n^{\frac{\phi_{\varphi}}{2}}),
\end{align}
since ${{{\sum_{i=1}^{\kappa_{\varphi}}d_{\varrho_{\phi,i}}}+2\(\upsilon_{\varphi}-\kappa_{\varphi}\)}}\leq 2{\phi_{\varphi}}$ and $\kappa_{\varphi}\geq 1.$
\end{itemize}
We have now established the following lemma.
\begin{lemma}\label{lemb1}
  For the $\varphi$-th connected graph ${\hat G_\varphi}$, if $\phi_\varphi>2$, then
  \begin{equation}
    \sum_{\hat G_\varphi}\E\prod_{t=1}^{\phi_\varphi}(\mb \bbA_{\hat G_{\varphi t}}\bs\varepsilon_{\hat G_{\varphi t}}-\E\mb \bbA_{\hat G_{\varphi t}}\bs\varepsilon_{\hat G_{\varphi t}})
    =o(n^{\phi_\varphi/2}),
    \label{eqlb1}
  \end{equation}
  and if $\phi_\varphi=2$, then
  \begin{equation}
    \sum_{\hat G_\varphi}\E\prod_{t=1}^{\phi_\varphi}(\mb \bbA_{\hat G_{\varphi t}}\bs\varepsilon_{\hat G_{\varphi t}}-\E\mb \bbA_{\hat G_{\varphi t}}\bs\varepsilon_{\hat G_{\varphi t}})
    =O(n^{\phi_\varphi/2})
    =O(n).
    \label{eqlb2}
  \end{equation}
\end{lemma}

Now, we return to the proof of the joint central limit theorem.

We have the following facts:
\begin{description}
  \item[(i) when $k$ is odd] Applying Lemma \ref{lemb1}, if $k$ is odd, since for any $G$ there are at least two connected subgraphs of $G$ that contain more than two constructing basic graphs, the second conclusion of (\ref{eqclt1}) holds.
  \item[(ii) when $k$ is even] When $k$ is even, $G$ consists of $u_{i,j}$ ($1\leq i, j\leq q$) connected subgraphs composed of two basic graphs $G_i$ if $i=j$ and one basic graph $G_i$ and one basic graph $G_j$ if $j\neq i$.
Clearly, we have $2u_{i,i}+\sum_{j\neq i} u_{i,j}=k_i$ and $u_{i,j}=u_{j,i}$.
Compare
\begin{equation}\label{eqclt3}
\sum_{G}\E\(
    \prod_{l=1}^q\prod_{\ell=1}^{k_l}(A_{G_{l,\ell}}\varepsilon_{{G_l,\ell}}-\E A_{G_{l,\ell}}\varepsilon_{{G_{l,\ell}}})\)
\end{equation}
with the expansion of
\begin{align}\label{eq218}
\prod_{1\leq i\leq j\leq q}\({\rm Cov(Q_i, Q_j)}\)^{u_{i,j}}=& (E(Q_1-EQ_1)^2)^{u_{1,1}}(E(Q_2-EQ_2)^2)^{u_{2,2}}\cdots (E(Q_q-EQ_q))^{u_{q,q}}\\\notag
                                                \times & (E(Q_1-EQ_1)E(Q_2-EQ_2))^{u_{1,2}}\cdots (E(Q_1-EQ_1)E(Q_q-EQ_q))^{u_{1,q}}\\\notag
                                                &\vdots \\\notag
                                                \times & (E(Q_{q-1}-EQ_{q-1})E(Q_q-EQ_q))^{u_{q-1,q}}.
\end{align}
 The latter expansion (\ref{eq218})
contains more terms than (\ref{eqclt3}), with more connections among the subgraphs. For the readers' convenience, in Figure \ref{difference}, we give an example where the term corresponding to $\hat G_{(1)}$ belongs to both (\ref{eq218}) and (\ref{eqclt3}) while the term corresponding to $\hat G_{(2)}$ belongs to (\ref{eq218}) but not to (\ref{eqclt3}).
\begin{figure}[H]
  \includegraphics[width=0.8\textwidth]{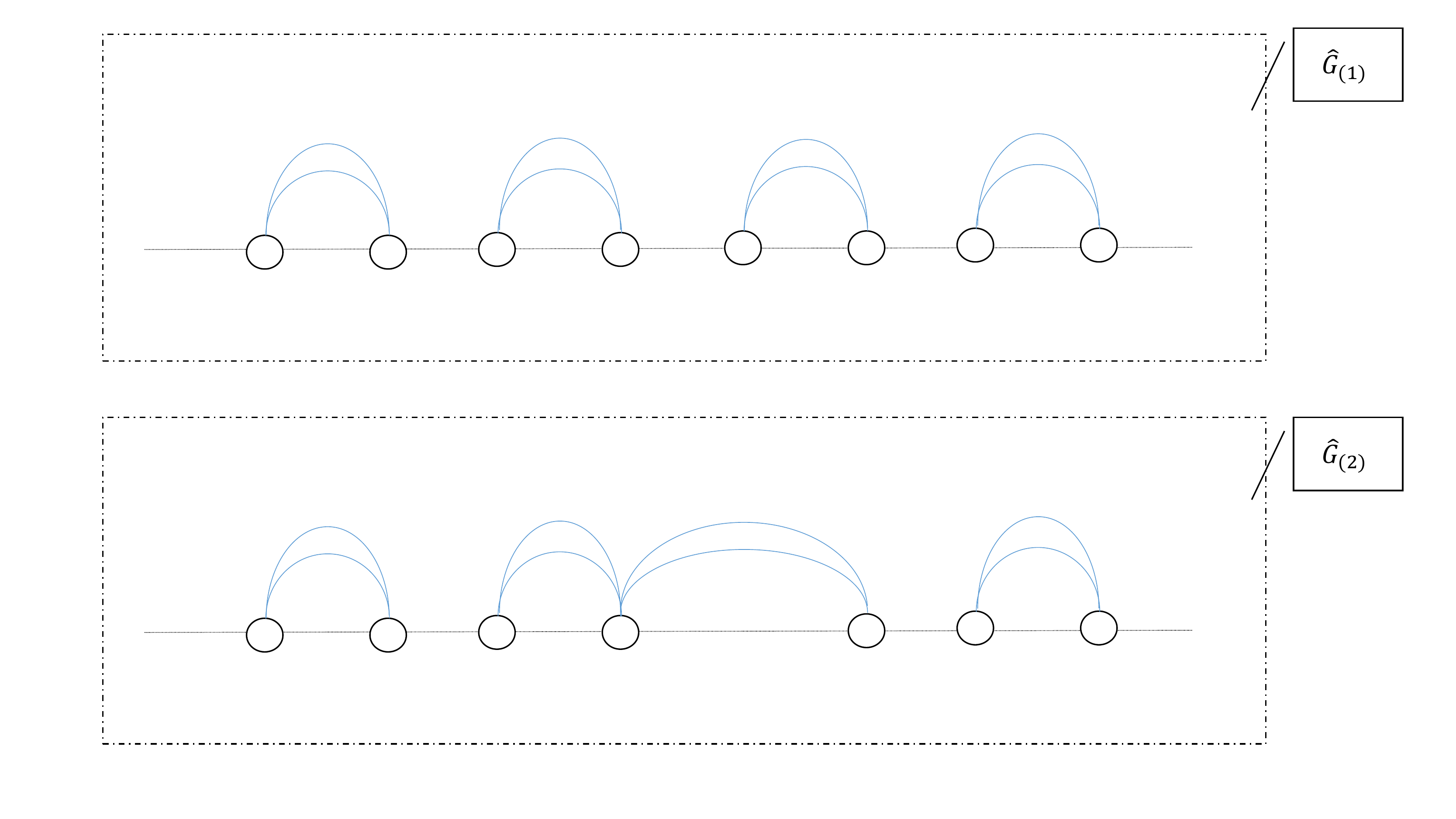}
  \caption{An example of ${\hat G_{(1)}}$ and ${\hat G_{(2)}}$.}\label{difference}	
\end{figure}

Therefore, for the same reason as the arguments above, we know that the whole contribution of the difference terms  between (\ref{eq218}) and (\ref{eqclt3})  has an order of $o(n^{k/2})$; thus,
\begin{equation}\label{eqclt4}
  \sum_{G}\E\(
    \prod_{l=1}^q\prod_{\ell=1}^{k_l}(A_{G_{l,\ell}}\varepsilon_{{G_l,\ell}}-\E A_{G_{l,\ell}}\varepsilon_{{G_{l,\ell}}})\)
  =\prod_{1\leq i\leq j\leq q}\({\rm Cov(Q_i, Q_j)}\)^{u_{i,j}}+o(n^{k/2}).
\end{equation}
Therefore, we obtain
\begin{align}
\E\prod_{l=1}^q(Q_l-\E Q_l)^{k_l}
=\sum_{\substack{2u_{1,1}+\sum_{j\neq 1} u_{1,j}=k_1\\ \vdots \\2u_{q,q}+\sum_{j\neq q} u_{q,j}=k_q}}\frac{k!}{2^{\sum_{l=1}^qu_{l,l}}\prod_{1\leq i\leq j\leq q}u_{i,j}!}\prod_{1\leq i\leq j\leq q}\({\rm Cov(Q_i, Q_j)}\)^{u_{i,j}}+o(n^{k/2}).
\label{eqclt5}
\end{align}

Substituting (\ref{eqclt5}) into (\ref{eqclt2}) yields
\begin{align}
&\E\(\sum_{l=1}^q\alpha_l \(Q_l-\E Q_l\)\)^k\\\notag
=&\sum_{\sum_{1\leq i\leq j\leq q}u_{i,j}=k}\frac{k!}{2^{\sum_{l=1}^qu_{l,l}}\prod_{1\leq i\leq j\leq q}u_{i,j}!}\prod_{l=1}^q\alpha_{l}^{2u_{l,l}+\sum_{j\neq l} u_{l,j}}\prod_{1\leq i\leq j\leq q}\({\rm Cov(Q_i, Q_j)}\)^{u_{i,j}}+o(n^{k/2})\\\notag
=&\frac{k!}{2^{k/2}(k/2)!}\(\sum_{i=1}^q\sum_{j=1}^q\alpha_{i}\alpha_{j}{\rm Cov}\(Q_{i},Q_{j}\)\)^{k/2}+o(n^{k/2})\\\notag
=&(k-1)!!\(\sum_{l_1=1}^q\sum_{l_2=1}^q\alpha_{l_1}\alpha_{l_2}{\rm Cov}\(Q_{l_1},Q_{l_2}\)\)^{k/2}+o(n^{k/2})\notag.\label{eqclt6}
\end{align}
\end{description}

Since we assume that there exists an $i$ such that ${\rm Cov}(Q_i,Q_i)$ has the same order as $n$, we conclude that for almost all $(\alpha_1,\alpha_2,\cdots,\alpha_q)\in \mathbb{R}^q$, we have that $$\(\sum_{l_1=1}^q\sum_{l_2=1}^q\alpha_{l_1}\alpha_{l_2}{\rm Cov}\(Q_{l_1},Q_{l_2}\)\)$$ has the same order as $n$. In fact, note that $n^{-1}\(\sum_{l_1=1}^q\sum_{l_2=1}^q\alpha_{l_1}\alpha_{l_2}{\rm Cov}\(Q_{l_1},Q_{l_2}\)\)$, and $f(\alpha_1,\alpha_2,\cdots,\alpha_q)$ is a polynomial in variables $\alpha_1,\alpha_2,\cdots,\alpha_q.$ We know from the fundamental properties of polynomials that one and exactly one of the following two cases holds: (1) The polynomials $f(\alpha_1,\alpha_2,\cdots,\alpha_q)\equiv 0$ for all vector $\(\alpha_1,\alpha_2,\cdots,\alpha_q\)$. (2) The Lebesgue measure of the set of vectors $\(\alpha_1,\alpha_2,\cdots,\alpha_q\)$ in the space $\mathbb{R}^q$ such that polynomials $f(\alpha_1,\alpha_2,\cdots,\alpha_q)=0$ is zero. We thus obtain the conclusion since (1) conflicts with our assumption by taking $\alpha_i=1$ and $\alpha_j=0$ for $j\neq i$.

Finally,  by applying the moment convergence theorem and continuity, we arrive at the fact that for all $\(\alpha_1,\alpha_2,\cdots,\alpha_q\)$, $$\sqrt{1/n}\(\sum_{l=1}^q\alpha_l\(Q_l-\E Q_l\)\)\stackrel{\mathcal D} \sim N(0,\frac{\(\sum_{l_1=1}^q\sum_{l_2=1}^q\alpha_{l_1}\alpha_{l_2}{\rm Cov}\(Q_{l_1},Q_{l_2}\)\)}{n}).$$

The proof of Theorem \ref{th:CLT} is complete.

\section{Conclusion and further discussion}
In this paper, we consider tests for detecting series correlation that are valid in both low- and high-dimensional linear regression models with random and fixed designs. The test statistics are based on the residuals of OLS and the residual maker matrix. We need no model assumptions on the regressor and/or dependent variable; thus, the tests are model-free. The asymptotic distribution of the statistics under the null hypothesis are obtained as a consequence of a general joint CLT of quadratic forms. Simulations are conducted to investigate the advantages of the proposed test procedures. The results show that the proposed tests perform well if $n-p$, where $n$ is the sample size and $p$ is the number of regressors, is not too small.

If we are concerned about the robustness, then we can use the standard residuals $\check\varepsilon_j = \varepsilon_j/\sqrt{p_{jj}}$ $1\le j\le
n$ instead of the original residuals.
 Then, the residual vector can be rewritten as $\check{\bs\varepsilon} =
(\check\varepsilon_1,\ldots,\check\varepsilon_n)'
= \mb D\mb R\mb \Sigma^{1/2} \bs\epsilon$, where
$\mb D$ is a diagonal matrix with diagonal entries $\{r_{jj}^{-1/2}\}_{j=1}^n$. Then, the test procedures remain valid after recalculating $m_{\tau}$ for $0\leq \tau\leq q$ and $v_{\tau_1\tau_2}$ for $0\leq \tau_1,\tau_2\leq q$ by replacing $\mb R$ with $\mb D\mb R$ in (\ref{mean}) and (\ref{cov}).

\end{document}